\documentclass[a4paper, 11pt]{amsart}
\usepackage{amssymb, enumitem, color, tikz, hyperref, cleveref}
\usepackage[T1]{fontenc}
\usepackage{mathptmx}
\usepackage{todonotes, comment}
\usepackage{array}
\usepackage{multirow}
\usepackage{graphicx}
\usepackage[margin=2.5cm]{geometry}

\newtheorem{theorem}{Theorem}
\numberwithin{theorem}{section}
\newtheorem{proposition}[theorem]{Proposition}
\newtheorem{lemma}[theorem]{Lemma}
\newtheorem{corollary}[theorem]{Corollary}

\theoremstyle{definition}
\newtheorem{definition}[theorem]{Definition}
\newtheorem{example}[theorem]{Example}
\newtheorem{remark}[theorem]{Remark}

\crefname{question}{question}{questions}

\numberwithin{equation}{section}

\DeclareMathOperator{\KK}{\mathbb{K}}
\DeclareMathOperator{\Ann}{\operatorname{Ann}}
\DeclareMathOperator{\Hom}{\operatorname{Hom}}

\DeclareMathOperator{\lcm}{\operatorname{lcm}}
\newcommand{\cR}{\mathcal R}
\newcommand{\Gr}{\operatorname{Gr}}
\newcommand{\Sym}{\operatorname{Sym}}
\newcommand{\rank}{\operatorname{rank}}

\newcommand{\curl}{\operatorname{curl}}

\title[Simultaneous apolarity and rational curves]{Simultaneous Waring problems for binary forms and the geometry of rational curves}

\author{Maciej Ga{\l}\k{a}zka \and Jaewoo Jung \and Hyunsuk Moon}

\address[Maciej Ga{\l}\k{a}zka]{Universit\"atsstraße 10, 78464 Konstanz, Germany}
\email{maciej.galazka@uni-konstanz.de}

\address[Jaewoo Jung]{DGIST BRL-AGSTA, 333 Techno Jungang-daero, Hyeonpung-eup, Dalseong-gun, Daegu, 42988, Republic of Korea}
\email{jaewoojung@dgist.ac.kr}

\address[Hyunsuk Moon]{School of Mathematics, Statistics and Data Science, Sungshin Women's University, Seoul 02844, Republic of Korea}
\email{hsmoon87@sungshin.ac.kr}

\thanks{All authors equally contributed to this work.}

\begin{document}

\begin{abstract}
We investigate the simultaneous Waring problem for binary forms and its geometric applications to rational curves.
By analyzing the minimal free resolution of the simultaneous annihilator, we establish formulas for both simultaneous Waring ranks and border Waring ranks including explicit computations for binary monomials.

As a geometric application, we explicitly describe parametrizations of rational curves up to projective equivalence.
We determine the minimum dimension of a rational normal scroll containing the curve as a section of its ruling in terms of the minimal generators of the corresponding simultaneous annihilator.
For rational curves lying as divisors on rational normal surface scrolls, we characterize their smoothness and describe the generators of their defining ideals.
\end{abstract}



\maketitle

\section{Introduction}
The polynomial Waring problem is a classical subject in algebraic geometry, seeking to determine the minimum number of linear forms required to express a given homogeneous polynomial of degree $d$ as a sum of their $d$-th powers.
This minimum is called the \emph{Waring rank} of the polynomial.
The \emph{border Waring rank} is defined using the Zariski closures of the loci of forms of bounded Waring rank.
Since the foundational works of Sylvester and Macaulay (cf. \cite{MR1735271} for the general setting and \cite{MR2754189} for the binary case), this problem has evolved into an active area of research, driven in part by its close connections to symmetric tensor decomposition (see, e.g., \cite{MR2447451}) and algebraic complexity theory (see, e.g., \cite{MR2628829,MR3729273}).

Throughout this paper, we work over an algebraically closed field $\mathbb{K}$ of characteristic zero.
To study the Waring decomposition of a homogeneous polynomial (or \emph{form}), the primary algebraic tool is its \emph{apolar ideal}.
Restricting our focus to binary forms, we identify $\Sym^d\mathbb{K}^2$ with $\mathcal{R}_d$, the degree-$d$ component of the polynomial ring $\mathcal{R}=\mathbb{K}[S,T]$.
By introducing the dual ring $R = \mathbb{K}[s,t]$, we can view $\mathcal{R}$ as an $R$-module via the differentiation action, and define the annihilator (or apolar ideal) of $F \in \mathcal{R}_d$ as 
\[
\Ann(F) = \{g \in R \mid g \circ F = 0\}.
\]
Sylvester classically showed that $\Ann(F)$ is generated by two forms $g_1$ and $g_2$ satisfying $\deg g_1 + \deg g_2 = d + 2$. 
Assuming $\deg g_1 \le \deg g_2$, the \emph{Waring rank} of $F$ is given by $\deg g_1$ if $g_1$ is square-free, and $\deg g_2$ otherwise, while its \emph{border Waring rank} is always equal to $\deg g_1$.

Geometrically, the Waring problem translates naturally to the study of secant varieties of Veronese varieties. 
In the binary setting, this picture reduces to secant varieties of the rational normal curve $\widetilde{\mathcal{C}}_d \subset \mathbb{P}_{\mathbb{K}}^d = \mathbb{P}(\mathcal{R}_d)$. 
A geometric problem arising in this context is to characterize the curve $\pi_{p}(\widetilde{\mathcal{C}}_d)$ obtained by a linear projection $\pi_p: \mathbb{P}_{\mathbb{K}}^d \dashrightarrow \mathbb{P}_{\mathbb{K}}^{d-1}$ from the center $p\in \mathbb{P}_{\mathbb{K}}^d$.
It turns out that the \emph{$\widetilde{\mathcal{C}}_d$-rank} of $p$, which is the smallest $q$ such that $q$-th secant variety of $\widetilde{\mathcal{C}}_d$ contains $p$, completely determines all the graded Betti numbers of the projected curve (cf. \cite{MR2299577}), and this rank exactly equals the minimal degree of the generators of $\Ann(F_p)$ where $F_p \in \mathcal{R}_d$ is the binary form corresponding to $p$. (See \cite[Section 2.2]{MR5047032} for details.)

A natural generalization in this context is the study of \emph{simultaneous Waring problem} (see, e.g., \cite{carlini_chipalkatti,angelini_galuppi_mella_ottaviani,carlini_ventura_simultaneous}).
Given homogeneous polynomials $F_1, \ldots, F_r$ with $\deg F_i = e_i$ for $i = 1, \dots, r$, one asks for the minimum number of linear forms that can simultaneously form a Waring decomposition for each $F_i$. In other words, we seek to determine
\begin{equation*}
  \min \left\{ k \in \mathbb{N} \;\middle|\; \exists \, l_1, \ldots, l_k \in \mathcal{R}_1 \text{ such that } F_i \in \langle l_1^{e_i}, \ldots, l_k^{e_i} \rangle \text{ for all } i = 1, \ldots, r \right\}.
\end{equation*}
This minimum is called the \emph{simultaneous Waring rank} of $F_1,\dots,F_r$.
The \emph{simultaneous border Waring rank} of $F_1,\dots,F_r$
is defined as
\begin{equation*}
\min\left\{k\in\mathbb N \;\middle|\; (F_1,\dots,F_r)\in \overline{ \bigcup_{l_1,\dots,l_k\in\mathcal R_1} \prod_{i=1}^r \langle l_1^{e_i},\dots,l_k^{e_i}\rangle} \right\},
\end{equation*}
where the bar denotes the Zariski closure in $\prod_{i=1}^r\mathcal R_{e_i}$.

Motivated by this, we aim to address the following questions in this article:
\begin{enumerate}
  \item Can we describe the algebraic structure of the simultaneous annihilator of $F_1, \ldots, F_r$?
  \item Can we give formulas for the simultaneous Waring rank and border Waring rank of $F_1, \ldots, F_r$ in terms of the generators of this annihilator?
  \item How can we understand the geometry of rational curves through the structure of the simultaneous annihilator?
\end{enumerate}

The first question is answered by using the classical structure theorems for graded rings of projective dimension $2$, namely the \emph{Hilbert-Burch theorem}.
As this foundational result serves as the starting point for our investigation, we highlight it as follows:
\begin{theorem}\label{thm:main}
Let $W \subset \mathcal{R}$ be a nonzero $R$-submodule minimally generated by the homogeneous polynomials $F_1, F_2, \dots, F_r \in \mathcal{R}$:
\[
W = R \circ F_1 + R \circ F_2 + \cdots + R \circ F_r.
\]
Then $I=\Ann(W)=\{g\in R~|~g\circ F=0 \text{~for all~} F\in W\}\subset R$ is an Artinian ideal minimally generated by $r+1$ binary forms $g_1,\ldots,g_{r+1}$, i.e.,
$$
\Ann(W) = (g_1, g_2, \dots, g_{r+1}).
$$
Conversely, let $I\subset R$ be an Artinian ideal minimally generated by the $r+1$ binary forms $g_1,\ldots,g_{r+1}$.
Then $W=I^{\perp}=\{F\in \mathcal{R}~|~g\circ F=0\text{~for all~} g\in I\}\subset \mathcal{R}$ is an $R$-submodule minimally generated by $r$ homogeneous polynomials $F_1,\ldots,F_r\in \mathcal{R}$.

Moreover, 
$$
\sum_{i = 1}^{r+1} \deg(g_i) = \sum_{j = 1}^{r} (\deg(F_j) + 2).
$$
\end{theorem}
   In this case, we always have $V(g_1,g_2,\ldots,g_{r+1}) = \{0\}$; see \Cref{lem:finitegeneration}.
The proof of \Cref{thm:main} is given in \Cref{sec:proofofMainThm}.

As an immediate consequence of \Cref{thm:main}, we explicitly characterize both the simultaneous and simultaneous border Waring ranks of binary forms in terms of the minimal generators of their annihilator ideal. 
\begin{theorem}\label{thm:simultaneousrank}
    Let $F_1, \dots, F_r \in \mathcal{R}$ be minimal homogeneous generators of the submodule $R \circ F_1 + \cdots + R \circ F_r$. 
    Let $g_1, \dots, g_{r+1} \in R$ be the minimal generators of $\Ann(F_1, \dots, F_r)$ ordered such that $d_1 \leq \dots \leq d_{r+1}$ where $d_i = \deg g_i$.
    Let $h_i = \gcd(g_1, \dots, g_i)$ for $1 \leq i \leq r+1$.
    Then the following hold:
    \begin{enumerate}
        \item \label{item:honestrank} The simultaneous Waring rank of $F_1, \dots, F_r$ is exactly $d_k$ where $k$ is the smallest integer such that the polynomial $h_k$ is square-free.
        \item \label{item:borderrank} The simultaneous border Waring rank of $F_1, \dots, F_r$ is exactly $d_1$.
    \end{enumerate}
\end{theorem}
The statements and proofs for \Cref{thm:simultaneousrank} are given in \Cref{cor:bordersimultaneousrank} and \Cref{thm:simultaneousrankalgorithm}.
Applying \Cref{thm:simultaneousrank}, we compute the simultaneous border ranks of binary monomials and simultaneous Waring ranks presented in \Cref{thm:monomialrank}.

The case where $F_1,\dots, F_r$ have the same degree is particularly well-studied due to the natural connection between the algebraic structure of simultaneous rank loci of symmetric tensors and the geometry of Grassmann secant varieties of Veronese varieties.
In this setting, the associated quotient algebra is known as a \emph{level algebra}, which has been actively investigated in the literature (e.g., \cite{MR3592466,MR2028071}).

With this motivations, we restrict our focus to binary forms of a fixed degree $d$ where the corresponding geometric counterpart is the rational normal curve $\widetilde{\mathcal{C}}_d$.
Any rational curve $\mathcal{C}$ parameterized by a birational map from $\mathbb{P}^1$ can be realized as a linear projection of $\widetilde{\mathcal{C}}_d$.
If we denote the center of this projection by $\mathbb{P}(L)$ for a subspace $L = \langle F_1,\dots, F_r \rangle \subset \mathbb{K}[S,T]_d$, then $\mathcal{C}$ is the image of $\widetilde{\mathcal{C}}_d$ under the rational map 
\[
\pi_{\mathbb{P}(L)}:\mathbb{P}(\mathrm{Sym}^d\mathbb{K}^2) \dashrightarrow \mathbb{P}(\mathrm{Sym}^d\mathbb{K}^2/L).
\]
Thus, it is natural to ask how the algebraic properties of $L$ dictate the geometry of $\mathcal{C}$.

To relate the structure of the annihilator and the parametrization of the curve, we consider an algebraic invariant $\rho(L)$ which we call the \emph{apolar defect}, defined as the number of minimal generators of $\Ann(L)$ in degree $d+1$.

The apolar defect $\rho(L)$ admits a natural algebraic interpretation as the dimension of the \emph{first prolongation} of $L$ (see \Cref{thm:inverse_space_top_generators}).
It also bounds the degree of the common factor $g$; that is, $\deg(g)\leq\rho(L)$ (see \Cref{cor:inequalityApolarDefect}).
As $L$ varies in the Grassmannian $\Gr(r,\mathcal{R}_d)$, we show that the locus $\mathfrak D_{\ge\tau}=\{L\in\Gr(r,\mathcal{R}_d)\mid\rho(L)\geq\tau\}$ is Zariski closed and locally defined by determinantal equations (see \Cref{prop:apolarDefectLocus}).

Using the apolar defect, we can describe the parametrization of the curve up to projective equivalence (see \Cref{thm:parametrizedCurve}).
Given $L = \langle F_1,\dots,F_r \rangle \subset \mathcal{R}_d$ with $\dim_{\mathbb{K}}L = r$ and $1\leq r \leq d-1$, suppose $g_1, \dots, g_{m}$ are the minimal generators of $\Ann(L)$ of degree at most $d$, where $m := r+1-\rho(L)$, and let $d_i := \deg(g_i)$.
Setting $g := \gcd(g_1, \dots, g_{m})$ and $g_i' := g_i/g$, the curve $\mathcal{C} = \pi_{\mathbb{P}(L)} (\widetilde{\mathcal{C}_d})$ is parametrized as follows:
\begin{equation}\label{equ:parametrizationCurve}
\left\{ \left[ \dots : g_i'(\alpha,\beta)\alpha^{d-d_i-j}\beta^{j} : \cdots \right]_{1 \le i \le m, \, 0 \le j \le d-d_i} \mid [\alpha:\beta] \in \mathbb{P}_{\mathbb{K}}^1 \right\} \subset \mathbb{P}^{d-r}_{\mathbb{K}}.
\end{equation}

An immediate consequence of this parametrization is that we can determine how many inner projections, counted with multiplicity, are involved in $\pi_{\mathbb{P}(L)}$ by computing the degree of the common factor $g=\gcd(g_1,\dots,g_m)$ (see \Cref{cor:detectInnerProjections}).
Furthermore, using the apolar defect $\rho(L)$, we construct a rational normal scroll of dimension $r+1-\rho(L)$ containing $\mathcal{C}$ as a section of its ruling and therefore improve the general upper bound $r+1$.
The bound $r+1$ follows from successive point projections using \cite[Theorem~7.3]{MR2274517} for the outer projections of the containing scrolls.

The following theorem establishes the minimality of this dimension under the section condition and provides a precise geometric interpretation of the apolar defect $\rho(L)$.
\begin{theorem}\label{thm:minimalSectionScroll}
Given the projected curve $\mathcal{C}=\pi_{\mathbb{P}(L)}(\widetilde{\mathcal{C}}_d)\subseteq\mathbb{P}_{\mathbb{K}}^{d-r}$ that is projectively equivalent to the curve with the parametrization in \eqref{equ:parametrizationCurve}, let $\varphi:\mathbb{P}_{\mathbb{K}}^1\longrightarrow\mathcal{C}$ be the morphism defined by \eqref{equ:parametrizationCurve}.
Let $m$ be an integer such that $2\leq m\leq \min\{r+1,d-r\}$.
Then $\rho(L)=r+1-m$ if and only if 
\begin{enumerate}
    \item $\varphi$ is birational onto its image and 
    \item $m$ is the smallest positive integer $n$ such that $\mathcal{C}$ is contained in an $n$-fold rational normal scroll in $\mathbb{P}_{\mathbb{K}}^{d-r}$ as a section of its ruling.
\end{enumerate}
\end{theorem}

We also explore more specialized geometric settings.
In particular, we study rational curves $\mathcal{C}$ that lie as divisors on a rational normal surface scroll $\mathcal{S}=S(a,b)$.
Here, $S(a,b)$ is obtained by joining corresponding points of two rational normal curves of degrees $a$ and $b$ in disjoint linear subspaces.
For a smooth scroll $\mathcal{S}=S(a,b)$ with $1\leq a\leq b$, we show that $\mathcal{C}$ is smooth whenever it is a section of the ruling of $\mathcal{S}$.
Conversely, if $\mathcal{C}$ is smooth and is not a rational normal curve, then it is a section of the ruling of $\mathcal{S}$ provided that $d$ is sufficiently large relative to $r$ (see \Cref{thm:surface-stable}).

Finally, as an application of our parametrization results, we study the defining ideal $I(\mathcal{C})$ of the curve $\mathcal{C}$. 
Assuming the morphism $\varphi$ is birational, \Cref{cor:RNC} shows that the maximal condition $\rho(L) = \dim(L)$ is equivalent to $\mathcal{C} = \overline{\pi_{\mathbb{P}(L)}(\widetilde{\mathcal{C}_d}\setminus \mathbb{P}(L))} \subset \mathbb{P}_{\mathbb{K}}^{d-r}$ being precisely the rational normal curve of degree $d-r$. 
When $\rho(L) = \dim_{\mathbb{K}}(L)-1$, we give precise conditions for the ideal $I(\mathcal{C})$ to be generated by quadrics.

\begin{theorem}\label{thm:definingequs}
Let $\mathcal{C}$ be the curve parameterized in \eqref{equ:parametrizationCurve} with $\rho(L) = \dim_{\mathbb{K}} (L) - 1$ (so that $m = 2$).
Let $d-r\geq 2$, $d_1\leq d_2$ and $d'=\deg(\gcd(g_1,g_2))$. Then the ideal $I(\mathcal{C})$ satisfies the following:
\begin{enumerate}
    \item\label{thm:definingequs_forms} $I(\mathcal{C})$ is generated by differences of products of linear forms.
    \item\label{thm:definingequs_quad} $I(\mathcal{C})$ is generated by quadratic forms if and only if 
    \[
    \begin{cases}
    d \ge 2r-d'+2 & \text{if } d_2 = d, \\
    d_1 \ge 2r-d'+2 & \text{otherwise}.
    \end{cases}
    \]
\end{enumerate}
\end{theorem}
We note that, if $d_2=d$ (equivalently, $d_1=r+1$), then $\mathcal{C}$ lies on a cone
over a rational normal curve of degree $d-r-1$.
If $d_2<d$, then $\mathcal{C}$ lies on a smooth rational normal surface scroll.

We provide explicit information about the numbers and degrees of minimal generators in \Cref{cor:classification}.
In the smooth-scroll case, the quadratic-generation criterion also follows from \cite[Theorem~4.4(2)]{MR3247023}.
Our approach gives a direct proof in terms of the minimal generators of the simultaneous annihilator, treating both smooth scrolls and cones.
Moreover, whenever $I(\mathcal C)$ is generated by quadrics,
it can be generated by quadrics of rank at most $4$.
Thus, \Cref{thm:definingequs} generalizes results on curves
of almost minimal degree (cf. \cite{MR5047032}).

In \Cref{sec:prel and main}, we give some preliminaries related to the apolar ideal of binary forms and the proof of \Cref{thm:main}.
For the second question, we provide an algorithm in \Cref{sec:simultaneous_rank} and exact rank formulas for binary monomials in terms of their exponents.
Regarding the third question, we provide general answers in \Cref{sec:rationalCurves}, and present further geometric and algebraic information regarding curves on surface scrolls in \Cref{sec:curvesinSurfacescroll}.

\section{Preliminaries and proof of \Cref{thm:main}}\label{sec:prel and main}
\subsection{Apolarity of binary forms}\label{subsec:apolarity}
In this subsection, we review some basics related to the notion of apolarity of binary forms. 
We fix some notation that will be maintained throughout this paper. 
 
Let $R:= \mathbb{K}[s,t]$ and $\mathcal{R}:= \mathbb{K}[S,T]$ be two bivariate standard graded polynomial rings over the field $\mathbb{K}$ with $\deg(s) = \deg(t) = \deg(S) = \deg(T) = 1$ so that $R = \bigoplus_{i\in \mathbb{Z}} R_i$ and $\mathcal{R} = \bigoplus_{i\in \mathbb{Z}} \mathcal{R}_i$ where $R_i := \bigoplus_{j = 0}^{i}\mathbb{K}s^jt^{i-j}$ and $\mathcal{R}_i := \bigoplus_{j = 0}^{i}\mathbb{K}S^jT^{i-j}$ denote the $i$-th graded components of $R$ and $\mathcal{R}$, respectively, for all $i\in \mathbb{Z}$. 
We call the polynomials $g\in R_i$ and $F \in \mathcal{R}_i$ binary forms of degree $i$ in $R$ and $\mathcal{R}$, respectively. 

\begin{definition}
For all $i,j\in \mathbb{Z}$, we consider the bilinear map of $\mathbb{K}$-vector spaces
\[
\circ_{(j,i)}:R_j\times \mathcal{R}_i \longrightarrow \mathcal{R}_{i-j} \quad\text{and the induced map}\quad \circ: R \times \mathcal{R} \longrightarrow \mathcal{R},
\]
where $\circ_{(j,i)}$ is defined by
\[
(g,F) \longmapsto g\circ_{(j,i)} F = g \circ F := \sum_{k=0}^j c_k\left(\frac{\partial}{\partial S}\right)^{j-k}\left(\frac{\partial}{\partial T}\right)^k F
\]
for $g=\sum_{k=0}^jc_k s^{j-k}t^k \in R_j$. 
The \textit{apolar ideal} (or annihilator) of a polynomial $F \in \mathcal{R}$ is defined to be the ideal consisting of all polynomials $g \in R$ for which $g \circ F = 0$, and is denoted by $\Ann(F)$.  
In other words, for $F\in\mathcal{R}$,
\[
\Ann(F) := \{g \in R : g \circ F = 0\}.
\]
\end{definition}
If the polynomial $F$ is homogeneous, then $\Ann(F)$ is a homogeneous ideal. Whenever an inverse-system module \(W\subset\mathcal R\) is regarded as a graded \(R\)-module, it is endowed with the inverse grading \(W_{-q}=W\cap\mathcal R_q\).

\begin{proposition} (cf. \cite[Theorem 1.44]{MR1735271})
    Suppose $F \in \mathcal{R}_d$ is a nonzero binary form of degree $d$. 
    Then, $\Ann(F)$ is a complete intersection ideal generated by two binary forms in $R$. 
    More explicitly,
    \[
    \Ann(F) = (g_1, g_2) \subset R
    \]
    for $g_1, g_2 \in R$ such that $V(g_1, g_2) = \emptyset$ in $\mathbb{P}_{\mathbb{K}}^1$ and $\deg g_1 + \deg g_2 = d + 2$. 
\end{proposition}

Apolarity is a powerful tool for studying the decomposition of forms into sums of powers of linear forms. 
This decomposition is classically known as a \emph{Waring decomposition} or a \emph{(generalized) additive decomposition} of forms. 
For a more comprehensive discussion on the generalized additive decomposition and apolarity, we refer the reader to \cite{MR1735271}.

\subsection{Duality of spaces via apolarity on binary forms}\label{subsec:spaceDuality}
For any integer $d \geq 0$, the pairing $\circ : R_d \times \mathcal{R}_d \to \mathcal{R}_0 = \mathbb{K}$ defines a non-degenerate inner product between the two vector spaces $R_d$ and $\mathcal{R}_d$. 
This means that for any basis $\{v_0, v_1, \dots, v_d\}$ of $R_d$ and any basis $\{w_0, w_1, \dots, w_d\}$ of $\mathcal{R}_d$, the matrix $v_{\bullet} \circ w_{\bullet} := (v_i \circ w_j \mid 0 \leq i,j \leq d) \in \mathbb{K}^{(d+1)\times(d+1)}$ has full rank $d+1$. 
Therefore, we can naturally identify $\mathcal{R}_d$ with the dual vector space of $R_d$ via this differential operation. 
In particular, for $F \in \mathcal{R}_d$ and $g \in R_d$, the inner product $g \circ F$ is explicitly given by
\[
g \circ F = \sum_{i = 0}^{d} c_i \left(\frac{\partial}{\partial S}\right)^{d-i} \left(\frac{\partial}{\partial T}\right)^{i} F(S,T) \in \mathbb{K},
\]
where $g = \sum_{i = 0}^{d} c_i s^{d-i} t^{i}$.

Let $\{s^d, s^{d-1}t, \dots, t^d\}$ be the standard monomial basis of $R_d$. 
The dual basis of $\mathcal{R}_d$ with respect to the inner product defined above is given by 
\[
\{S^{[d]}, \dots, S^{[d-i]}T^{[i]}, \dots, T^{[d]}\} := \{\frac{S^d}{d!}, \dots, \frac{S^{d-i}}{(d-i)!}\frac{T^i}{i!}, \dots, \frac{T^d}{d!}\}. 
\]
Indeed, we have
\begin{equation*}\label{eqn:duality}
s^{d-i}t^i\circ S^{[d-j]}T^{[j]}=\delta_{ij},
\end{equation*}
where $\delta_{ij}$ denotes the Kronecker delta.

\begin{lemma}(\cite[Lemma 2.3]{MR5047032})\label{lem:evaluation}
Let $d$ be a positive integer.
For any binary forms $g(s,t)\in \mathbb{K}[s,t]_d$ and $(\alpha S+\beta T)^d\in \mathbb{K}[S,T]_d$, we have
\[
g(s,t)\circ (\alpha S+\beta T)^d = d!g(\alpha,\beta).
\]
\end{lemma}
\Cref{lem:evaluation} implies that the inner product of $g(s,t) \in \mathbb{K}[s,t]_d$ with $(\alpha S + \beta T)^d$ corresponds to evaluating the form $g(s,t)$ at the point $(\alpha, \beta)$.

We identify the projective space $\mathbb{P}^d$ with the projectivization $\mathbb{P}(\mathbb{K}[S,T]_d)$ of the space of binary forms of degree $d$.
More explicitly, a point $p = [p_0 : p_1 : \dots : p_d] \in \mathbb{P}^d$ corresponds to the binary form (up to a scalar factor)
\[
F_p = \sum_{i=0}^d p_i S^{[d-i]} T^{[i]} = \frac{1}{d!} \sum_{i=0}^d p_i \binom{d}{i} S^{d-i} T^i \in \mathbb{K}[S,T]_d,
\]
where $(p_0, \dots, p_d)$ are the coefficients of $F_p$ with respect to the dual basis $\{S^{[d-i]}T^{[i]} \mid 0 \leq i \leq d\}$ of $\mathcal{R}_d$.
Under this identification, the torus-fixed point $e_i = [\delta_{i0} : \dots : \delta_{id}] \in \mathbb{P}^d$ corresponds to $S^{[d-i]}T^{[i]}$, and a point $\nu_d([\alpha:\beta])$ on the standard rational normal curve $\widetilde{\mathcal{C}}_d = \nu_d(\mathbb{P}^1) \subset \mathbb{P}^d$ corresponds to the pure power $(\alpha S + \beta T)^d$ of linear form (up to a nonzero scalar).

Under the apolar pairing $\circ : R_d \times \mathcal{R}_d \to \mathbb{K}$, the orthogonal complement of $F_p$ is precisely its degree-$d$ annihilator, i.e., $(F_p)^\perp = (\Ann(F_p))_d$. 
More generally, for any linear subspace $L = \langle F_1, \dots, F_r \rangle \subseteq \mathbb{K}[S,T]_d$, we have
\[
L^\perp = \bigcap_{i=1}^r F_i^\perp = \Ann(L)_d.
\]

\subsection{General facts on the apolarity of forms}\label{subsec:general_apolarity}
In this subsection, we review and establish several fundamental properties regarding the apolar ideal and the associated $R$-submodule $W$.
Throughout this subsection, we let $R$ be a standard graded polynomial ring over $\mathbb{K}$ and let $\mathcal{R}$ be its dual polynomial ring on which $R$ acts via differentiation.
The results in this subsection hold for general graded polynomial rings in an arbitrary number of variables.

\begin{lemma}\label{lem:ann_equality}
Let $F_1, \dots, F_r \in \mathcal{R}$ be homogeneous polynomials, and let $W = R \circ F_1 + \dots + R \circ F_r$ be the $R$-submodule generated by $F_1, \dots, F_r$. 
Then, we have 
\[
\Ann(W) = \Ann(F_1, \dots, F_r) = \bigcap_{i=1}^r \Ann(F_i),
\]
where $\Ann(W) := \{g \in R \mid g \circ F = 0 \text{ for all } F \in W\}$ and $\Ann(F_1, \dots, F_r) := \{g \in R \mid g \circ F_i = 0 \text{ for } i = 1, \dots, r\}$.
\end{lemma}
\begin{proof}
The inclusion $\Ann(W) \subseteq \Ann(F_1, \dots, F_r) $ is immediate from the definition. 
For the converse, suppose $g \in \Ann(F_1, \dots, F_r)$. 
That is, $g \circ F_i = 0$ for each $i = 1, \dots, r$. 
Since the differential operators defining the $R$-action commute, for any element $F = \sum_{i=1}^r h_i \circ F_i \in W$ with $h_i \in R$, we have 
\[
g \circ (h_i \circ F_i) = h_i \circ (g \circ F_i) = 0.
\]
Thus, $g \circ F = 0$ for all $F \in W$ which implies $g \in \Ann(W)$.
\end{proof}

\begin{lemma}\label{lem:finitegeneration}
Let $W$ be the nonzero graded $R$-submodule in $\mathcal{R}$ and let $\Ann(W) \subset R$ be its apolar ideal. 
The following conditions are equivalent:
\begin{enumerate}
    \item\label{lem:finitegeneration1} $V(\Ann(W)) = \{0\} \subset \mathbb{A}^n$.
    \item\label{lem:finitegeneration2} There exists an integer $d \in \mathbb{N}$ such that $\Ann(W)_d = R_d$.
    \item\label{lem:finitegeneration3} $W$ is a finitely generated $R$-submodule.
\end{enumerate}
\end{lemma}

\begin{proof}
The equivalence \eqref{lem:finitegeneration1} $\iff$ \eqref{lem:finitegeneration2} follows immediately from the (homogeneous) Nullstellensatz.

\medskip
\noindent \eqref{lem:finitegeneration3} $\implies$ \eqref{lem:finitegeneration2}: Suppose that $W$ is finitely generated. 
Let $\{F_1, \dots, F_r\}$ be a minimal set of homogeneous generators for $W$, and choose an integer $d$ such that $d > \max\{\deg(F_i) \mid i = 1, \dots, r\}$. 
For any form $g \in R_d$, we have $g \circ F_i = 0$ for each $i = 1, \dots, r$ since $\deg(g) = d > \deg(F_i)$. 
This implies that $g \circ F = 0$ for all $F \in W$, yielding $g \in \Ann(W)_d$. 
Thus, we have $R_d = \Ann(W)_d$.

\medskip
\noindent \eqref{lem:finitegeneration2} $\implies$ \eqref{lem:finitegeneration3}: 
Assume that $\Ann(W)_d = R_d$ for some $d \in \mathbb{N}$. 
This condition implies that $R_{\geq d} \subseteq \Ann(W)$, and thus the quotient ring $R/\Ann(W)$ is an Artinian ring. 
Under Macaulay's inverse system correspondence, the $R$-submodule $W$ is isomorphic to the graded $\mathbb{K}$-dual of $R/\Ann(W)$. 
Since $R/\Ann(W)$ is finite-dimensional $\mathbb{K}$-vector space, $W$ is finitely generated as an $R$-module.
\end{proof}

\begin{lemma}\label{lem:minimality_cond}
Let $F_1,\dots,F_r \in \mathcal{R}$ be forms generating the $R$-module $W$.
Then $\{F_1, \dots, F_r\} \subset \mathcal{R}$ is a minimal generating set if and only if, for each $i = 1, \dots, r$,
\begin{equation}\label{equ:minimality_cond}
    \bigcap_{j \neq i} \Ann(F_j) \nsubseteq \Ann(F_i).
\end{equation}
In particular, if $F_1, \dots, F_r \in \mathcal{R}_d$ are forms of the same degree $d$, this condition is satisfied if and only if the set $\{F_1, \dots, F_r\}$ is linearly independent in $\mathcal{R}_d$.
\end{lemma}
\begin{proof}
Suppose that $\{F_1, \dots, F_r\}$ is a minimal generating set.
For the sake of contradiction, suppose that $\bigcap_{j \neq i} \Ann(F_j) \subseteq \Ann(F_i)$ for some $i$. 
This inclusion implies that $\bigcap_{j \neq i} \Ann(F_j) = \bigcap_{j=1}^r \Ann(F_j)$. This contradicts the assumption that $\{F_1, \dots, F_r\}$ is a minimal generating set.

Conversely, suppose that \eqref{equ:minimality_cond} holds for every $i$.
If the generating set is not minimal, then $F_i = \sum_{j \neq i} h_j \circ F_j$ for some $i$ and some $h_j \in R$.
For any $g \in \bigcap_{j \neq i} \Ann(F_j)$, the commutativity of differential operators gives $g \circ F_i = \sum_{j \neq i} h_j \circ (g \circ F_j) = 0$.
Thus $\bigcap_{j \neq i} \Ann(F_j) \subseteq \Ann(F_i)$ which is a contradiction.

Finally, if all $F_i$ have degree $d$, the degree-$d$ part of $\sum_{j \neq i} R \circ F_j$ is precisely the $\mathbb{K}$-linear span of the other generators.
Hence minimality is equivalent to linear independence.
\end{proof}

In view of \Cref{lem:minimality_cond}, the condition \eqref{equ:minimality_cond} can be utilized as an alternative characterization of the minimality of a generating set for the $R$-submodule $W$.

\subsection{Canonical modules and dual resolutions}\label{subsec:canonical_modules}
We recall the necessary algebraic foundations regarding canonical modules and dual resolutions from \cite{MR1322960}.

\begin{definition}
Let $A$ be a finite-dimensional graded $\mathbb{K}$-algebra with $A_0=\mathbb{K}$.
The \emph{graded canonical module} of $A$, denoted by $\omega_A$, is defined as the graded $\mathbb{K}$-dual $\operatorname{Hom}_{\mathbb{K}}(A,\mathbb{K})$ where $(\omega_A)_j=\operatorname{Hom}_{\mathbb{K}}(A_{-j},\mathbb{K})$.
\end{definition}

\begin{definition}[cf. {\cite[Section~21.11]{MR1322960}}]
Let $A$ be a positively graded Cohen--Macaulay ring
with $A_0=\mathbb{K}$.
A finitely generated graded $A$-module $W$ is said to be a \emph{graded canonical module} for $A$ if there exists a homogeneous element $x\in A$ of positive degree $\delta$ that is a non-zerodivisor on both $A$ and $W$ such that $(W/xW)[\delta]$ is a graded canonical module for $A/xA$.
\end{definition}

We recall the graded version of \cite[Theorem 21.15]{MR1322960}:
\begin{theorem}\label{thm:eisenbud_ext}
Let $R$ be a graded polynomial ring with canonical module $\omega_R = R[-\dim R]$. If $A = R/I$ for some homogeneous ideal $I$, and $A$ is Cohen-Macaulay, then $A$ admits a canonical module. In fact, if $c = \dim R - \dim A$, then
\[
\omega_A = \operatorname{Ext}^c_R(A, \omega_R).
\]
\end{theorem}

\begin{example}
The canonical module for the standard bivariate polynomial ring $R = \mathbb{K}[s,t]$ is $\omega_R = R[-2]$.
\end{example}

We recall the graded version of \cite[Corollary 21.16]{MR1322960}:
\begin{corollary}\label{corollary:dual_resolutions}
Let $R$ be a graded polynomial ring and let $I$ be a homogeneous ideal of codimension $c$ in $R$. 
Suppose that $A = R/I$ is Cohen-Macaulay. 
If $\mathcal{F}$ is a minimal graded free resolution of $A$ as an $R$-module, then the length of $\mathcal{F}$ is $c$ and its dual complex $\mathcal{F}^* = \operatorname{Hom}_R(\mathcal{F}, \omega_R)$ is a minimal graded free resolution of the canonical module $\omega_A$.
\end{corollary}

\begin{theorem}\label{thm:canonical_w}
Let $W$ be the $R$-submodule of $\mathcal{R}$ minimally generated by forms $F_1, \dots, F_r$.
Then $W$ is the canonical module of $R/I$, where $I = \bigcap_{i=1}^r \Ann(F_i)$. In particular, $W$ and $R/I$ possess dual minimal graded free resolutions.
\end{theorem}

\begin{proof}
By \Cref{lem:ann_equality}, we have $I = \Ann(W)$, and from \cite[Lemma A.4]{MR4634294}, it follows that $W \cong \operatorname{Hom}_{\mathbb{K}}(R/I, \mathbb{K})$. Since $R/I$ is zero-dimensional, this graded $\mathbb{K}$-dual coincides with the canonical module $\omega_{R/I}$. The assertion regarding dual resolutions follows immediately from \Cref{corollary:dual_resolutions}.

\end{proof}

\subsection{Structures of simultaneous annihilator of binary forms}\label{sec:proofofMainThm}
In this subsection, we will give a proof of \Cref{thm:main} and an example illustrating how the proof of \Cref{thm:main} works.

\begin{proof}[Proof of \Cref{thm:main}]
Note that $\Ann(W)$ is a finitely generated homogeneous ideal in $R = \mathbb{K}[s,t]$. 
Let $\Ann(W) = (g_1, \dots, g_m)$ be a minimal generating set for some $m \in \mathbb{N}$ and $d_i = \deg(g_i)$ for $i = 1, \dots, m$. 
By the Hilbert-Burch theorem (cf. \cite[Proposition 3.8]{MR2103875}), since $R/\Ann(W)$ is a zero-dimensional Cohen-Macaulay ring, the minimal graded free resolution of $R/\Ann(W)$ is of the form
\[
0 \longrightarrow \bigoplus_{j = 1}^{m-1} R[-e_j] \longrightarrow \bigoplus_{i = 1}^{m} R[-d_i] \longrightarrow R \longrightarrow R/\Ann(W) \longrightarrow 0,
\]
where the shifts $\{-d_i\}$ and $\{-e_j\}$ satisfy the relation $\sum_{i = 1}^{m} d_i = \sum_{j = 1}^{m-1} e_j$.
Applying $\operatorname{Hom}_R(-,\omega_R)$ to this free resolution, where $\omega_R = R[-2]$, gives a dual complex. 
Since we have 
\[
\operatorname{Ext}_R^2(R/\Ann(W),R[-2]) \cong \omega_{R/\Ann(W)},
\]
\Cref{thm:eisenbud_ext} and \Cref{corollary:dual_resolutions} give the minimal free resolution
\[
0 \longrightarrow R[-2] \longrightarrow \bigoplus_{i = 1}^{m} R[d_i-2] \longrightarrow \bigoplus_{j = 1}^{m-1} R[e_j-2] \longrightarrow \omega_{R/\Ann(W)} \longrightarrow 0.
\] 
Furthermore, by \Cref{thm:canonical_w}, we identify $\omega_{R/\Ann(W)} \cong W$.

Since $W$ is minimally generated by exactly $r$ elements $\{F_1, \dots, F_r\}$, the rank of the free module corresponding to the generators of $W$ must satisfy $m - 1 = r$. 
Thus, $\Ann(W)$ is minimally generated by $r+1$ elements $g_1, \dots, g_{r+1}$. 
Moreover, since $\deg(F_j) = e_j - 2$ for $j = 1, \dots, r$,
\[
\sum_{i = 1}^{r+1} \deg(g_i) = \sum_{j = 1}^{r} (\deg(F_j) + 2).
\]
The converse follows by applying the preceding argument to the Macaulay inverse system $W=I^\perp$ since $\Ann(W)=I$ by Macaulay inverse-system duality.
\end{proof}
\begin{example}
   Let $F_1 =  S^{[5]} T^{[3]}$ and $F_2 = S^{[2]} T^{[4]}$. Then
   \begin{equation*}
     \Ann(F_1,F_2) = (s^6, s^3 t^4, t^5).
   \end{equation*}
   Set $W = R\circ F_1 + R\circ F_2$.
   The resolution of $R/\Ann(W)$ is given as follows:
   \begin{equation*}
     0 \to R[-10]\oplus R[-8] \xrightarrow{M} R[-6]\oplus R[-7] \oplus R[-5] \to R \to R/\Ann(W) \to 0,
   \end{equation*}
   where
   \begin{equation*}
       M = \begin{bmatrix}
           t^4 & 0 \\
           -s^3 & t \\
           0 & -s^3
       \end{bmatrix}.
   \end{equation*}
   By applying the $\Hom_R(-, R[-2])$ functor, we obtain the dual resolution: 
   \begin{equation*}
     0 \to R[-2]\to R[4]\oplus R[5] \oplus R[3] \xrightarrow{M^T} R[8]\oplus R[6] \to W \to 0.
   \end{equation*}
   This is the minimal resolution of the module $W = R\circ F_1 + R\circ F_2$.
\end{example}
For a general $R$-submodule $W$ generated by binary monomials, the explicit description of the minimal free resolution of $\Ann(W)$ and its dual is provided following \Cref{lem:monomialsimultaneousrank}.

\section{Simultaneous Waring rank of forms}\label{sec:simultaneous_rank}
As an immediate application of \Cref{thm:main}, this section investigates the simultaneous Waring rank and border rank of a collection of binary forms. 
We first characterize the simultaneous border rank using the lowest degree generator of the joint annihilator ideal.
\begin{corollary}\label{cor:bordersimultaneousrank}
    Let $F_1, \dots, F_r\in \mathcal{R}$ be homogeneous minimal generators of $R\circ F_1 + \cdots + R\circ F_r$ with degrees $e_1,\ldots,e_r$ and $g$ be a non-zero generator of the lowest degree in $\operatorname{Ann}(F_1, \dots, F_r)$. 
    Then the simultaneous border Waring rank of $F_1, \dots, F_r$ is equal to $\deg g$.
    
    In particular, if a lowest-degree generator $g$ of $\Ann(F_1,\dots,F_r)$ is square-free, then the simultaneous Waring rank of $F_1, \dots, F_r$ also equals $\deg g$.
\end{corollary}

\begin{proof}
    Let $\deg g = d$. 
    Since $g$ is a binary form over $\KK$, we can factor it into linear forms as $g = \prod_{j = 1}^{q} \ell_j^{a_j} \in \Ann(F_1, \dots, F_r)$ where $\sum_{j=1}^q a_j = d$. 
    By the (generalized) apolarity (cf. \cite[Proposition 1.36]{MR1735271}), there exists a unique generalized additive decomposition of $F_i$ for each $i = 1, \dots, r$:
    \[
    F_i = \sum_{j = 1}^{q} l_j^{e_i-a_j+1} G_{ij},
    \]
    where $l_j\in \mathcal{R}_1$ is the annihilator of $\ell_j$ and $G_{ij} \in \mathcal{R}_{a_j - 1}$. 
    This factorization directly implies that the simultaneous border Waring rank of $F_1, \dots, F_r$ is equal to $d$ since there is no form in $\Ann(F_1,\ldots,F_r)$ with degree less than $d$.
    
    Moreover, if $g$ has only simple roots (i.e., $a_j = 1$ for all $j = 1, \dots, q$), then we obtain the standard simultaneous decomposition for each $F_i$:
    \[
    F_i = \sum_{j = 1}^{d} c_{ij} l_j^{e_i}.
    \]
    Hence, the simultaneous Waring rank of $F_1, \dots, F_r$ also equals $d$ since the simultaneous Waring rank is greater than or equal to the simultaneous border Waring rank.
\end{proof}

Next, we provide an algorithmic criterion to compute the exact simultaneous Waring rank. 
We begin by introducing a lemma necessary for its proof.
\begin{lemma}\label{lem:honestsimultaneousrank}
    Let $I = (g_1, \dots, g_r)$ be a homogeneous ideal in $R = \mathbb{K}[s,t]$ and let $d_i = \deg g_i$ for $1 \leq i \leq r$. 
    For any integer $d \geq \max\{d_1, \dots, d_r\}$, let $\{h_1, \dots, h_q\}$ be a basis for the vector space $I_d$.
    \begin{enumerate}
        \item\label{lem:same_gcd} Then the greatest common divisors satisfy
        \[
        \gcd(g_1, \dots, g_r) = \gcd(h_1, \dots, h_q)
        \]
        up to a nonzero scalar.
        \item\label{lem:bertini} Suppose $\gcd(h_1, \dots, h_q)$ is square-free. 
        Then a general element of the linear space $\langle h_1, \dots, h_q \rangle = I_d$ is square-free.
    \end{enumerate}
\end{lemma}

\begin{proof}
\eqref{lem:same_gcd} The assertion follows from the fact that the graded component $I_d$ is spanned by the elements of the form $R_{d-d_i} \cdot g_i$ for $1 \leq i \leq r$. 
Therefore,
\[
\gcd(h_1, \dots, h_q) = \gcd(R_{d-d_1}\cdot g_1, \dots, R_{d-d_r}\cdot g_r) = \gcd(g_1, \dots, g_r).
\]
    
\eqref{lem:bertini} Let $h = \gcd(h_1, \dots, h_q)$ and $h_i' = h_i/h$ for each $i = 1, \dots, q$. 
By construction, the linear system associated with $\langle h_1', \dots, h_q'\rangle$ is base-point-free. 
Consequently, by Bertini's theorem (see \cite[Chapter III Corollary~10.9]{hartshorne}), the vanishing scheme $V(p)$ of a general element $p \in \langle h_1', \dots, h_q'\rangle$ is reduced, meaning that $p$ is square-free.
Since $h$ is square-free and shares no common roots with a general $p$, a general element $h \cdot p \in \langle h_1, \dots, h_q \rangle$ is also square-free.
\end{proof}

We now present an algorithmic theorem to determine the simultaneous Waring rank of binary forms.
\begin{theorem}\label{thm:simultaneousrankalgorithm}
    Let $F_1, \dots, F_r\in \mathcal{R}$ be homogeneous minimal generators of $R\circ F_1 + \cdots + R\circ F_r$ and $g_1, \dots, g_{r+1} \in R$ be the minimal generators of $\Ann(F_1, \dots, F_r)$ ordered such that $\deg g_1 \leq \deg g_2 \leq \dots \leq \deg g_{r+1}$.
    For each $1 \leq i \leq r+1$, let $h_i = \gcd(g_1, \dots, g_i)$. 
    Then the simultaneous Waring rank of $F_1, \dots, F_r$ is equal to $\deg g_k$ where $k$ is the smallest integer such that $h_k$ is square-free.
\end{theorem}
\begin{proof}
    Let $\deg g_k = d$.
    If $h_k$ is square-free, then a general polynomial in the graded component $(\Ann(F_1,\dots,F_r))_d$ is square-free by \Cref{lem:honestsimultaneousrank}\eqref{lem:bertini}.
    Thus, any such polynomial gives a simultaneous Waring rank decomposition of length $d$ by the Apolarity Lemma (cf. \cite[Lemma 1.31]{MR1735271}).

    On the other hand, if $e<\deg g_1$, then $\Ann(F_1,\dots,F_r)_e=0$, so no simultaneous Waring decomposition of length $e$ exists. 
    For $\deg g_1\leq e < d$, let $j$ be the largest integer such that $\deg g_j\leq e$.
    By the minimality of $k$, $h_j$ is not square-free. 
    Since $h_j$ is the greatest common divisor of the entire graded component $\Ann(F_1, \dots, F_r)_e$ according to \Cref{lem:honestsimultaneousrank}\eqref{lem:same_gcd}, every form of degree $e$ in $\Ann(F_1, \dots, F_r)$ is necessarily divisible by a square. 
    Hence, by the Apolarity Lemma, there does not exist a simultaneous Waring decomposition of length $e$.
\end{proof}
\begin{remark}\label{rem:simultaneousCatalecticants}
The simultaneous Waring rank and a corresponding decomposition can also be determined directly using \emph{stacked catalecticant matrices}.
In each degree $d$, the single annihilator $\Ann(F_i)_d$ is obtained by calculating the kernel of the catalecticant map
$R_d \to \mathcal{R}_{e_i-d}$.
Therefore, $\Ann(F_1,\ldots,F_r)_d$ is the common kernel of these maps which can be found by stacking the corresponding matrices vertically and taking the kernel of the resulting block matrix.
\end{remark}
As an application of \Cref{cor:bordersimultaneousrank} and the algorithmic result in \Cref{thm:simultaneousrankalgorithm}, we determine both the simultaneous rank and simultaneous border rank of binary monomials.
To apply these results, we first characterize the conditions on the exponents of binary monomials for them to minimally generate an $R$-submodule $W$, and subsequently determine the minimal generators of its annihilator.
\begin{lemma}\label{lem:monomialsimultaneousrank}
Let $F_i := S^{[a_i]}T^{[b_i]}$ for $i = 1, \dots, r$ be binary monomials ordered such that $a_1 \leq \dots \leq a_r$. 
\begin{enumerate}
    \item\label{lem:monomialMinGensCond} The $R$-module $W := R \circ F_1 + \dots + R \circ F_r$ is minimally generated by $\{F_1, \dots, F_r\}$ if and only if $a_1 < \dots < a_r$ and $b_1 > \dots > b_r$.
    \item\label{lem:monomialAnn} If $W$ is minimally generated, then
    \[
    \Ann(F_1, \dots, F_r) = (t^{b_1+1}, s^{a_1+1}t^{b_2+1}, \dots, s^{a_{r-1}+1}t^{b_r+1}, s^{a_r+1}).
    \]
\end{enumerate}
\end{lemma}
\begin{proof} 
\eqref{lem:monomialMinGensCond} Suppose that $W$ is minimally generated by $\{F_1, \dots, F_r\}$. 
First, assume for contradiction that $a_k = a_{k+1}$ for some $1 \leq k \leq r-1$. 
Without loss of generality, we may assume $b_k \leq b_{k+1}$. 
The assumption immediately implies that $\Ann(F_{k+1}) \subset \Ann(F_k)$. 
Thus, we have
\[
\bigcap_{i \neq k} \Ann(F_i) \subset \Ann(F_{k+1}) \subset \Ann(F_k),
\]
which contradicts the minimality of the generating set by \Cref{lem:minimality_cond}. 
Thus, we must have $a_1 < \dots < a_r$.

Next, assume for contradiction that $b_k \leq b_{k+1}$ for some $1 \leq k \leq r-1$. 
Since $a_k < a_{k+1}$, we have $\Ann(F_{k+1}) \subset \Ann(F_k)$.
This implies $\bigcap_{i \neq k} \Ann(F_i) \subset \Ann(F_k)$ which contradicts the minimality of the generators. 
Therefore, we conclude that $b_1 > \dots > b_r$.

Conversely, suppose that $a_1 < \dots < a_r$ and $b_1 > \dots > b_r$. 
By \Cref{lem:minimality_cond}, it suffices to show that $\bigcap_{j \neq i} \Ann(F_j) \not\subseteq \Ann(F_i)$ for each $i = 1, \dots, r$. 
For a fixed index $i$, we claim that $g_i := s^{a_i}t^{b_i}$ belongs to $\bigcap_{j \neq i} \Ann(F_j)$ but not to $\Ann(F_i)$. 
Indeed, it is clear that $g_i \notin \Ann(F_i)$ since $g_i \circ F_i = 1 \neq 0$. 
On the other hand, let $j \neq i$ be any other index.
\begin{itemize}
    \item If $j < i$, then $g_i \circ F_j = s^{a_i}t^{b_i} \circ S^{[a_j]}T^{[b_j]} = s^{a_i} \circ S^{[a_j]}T^{[b_j-b_i]} = 0$.
    \item If $j > i$, then $g_i \circ F_j = s^{a_i}t^{b_i} \circ S^{[a_j]}T^{[b_j]} = t^{b_i} \circ S^{[a_j-a_i]}T^{[b_j]} = 0$.
\end{itemize}
Thus, $g_i \in \bigcap_{j \neq i} \Ann(F_j) \setminus \Ann(F_i)$ for all $i$.

\eqref{lem:monomialAnn} We proceed by induction on $r$. 
For $r=1$, we have $\Ann(F_1) = (t^{b_1+1}, s^{a_1+1})$.
Now, assume that the statement holds for $r=k \geq 1$.
That is,
\[
\Ann(F_1,\dots,F_k) = (t^{b_1+1}, s^{a_1+1}t^{b_2+1}, \dots, s^{a_{k-1}+1}t^{b_k+1}, s^{a_k+1}).
\]
For convenience, let $a_0 = -1$. 
By \Cref{lem:ann_equality}, we have
\[
\Ann(F_1,\dots, F_k ,F_{k+1}) = \Ann(F_1,\dots,F_k) \cap \Ann(F_{k+1})
\]
Since the intersection of monomial ideals is generated by the least common multiples (lcm) of all pairs of their generators (see, e.g., \cite[Proposition 1.2.1]{MR2724673}), $\Ann(F_1,\dots,F_{k+1})$ is generated by the following elements:
\begin{align*}
    \lcm(s^{a_{j-1}+1}t^{b_j+1}, s^{a_{k+1}+1}) &= s^{a_{k+1}+1}t^{b_j+1} \quad \text{for } j = 1, \dots, k, \\
    \lcm(s^{a_k+1}, s^{a_{k+1}+1}) &= s^{a_{k+1}+1}, \\
    \lcm(s^{a_{j-1}+1}t^{b_j+1}, t^{b_{k+1}+1}) &= s^{a_{j-1}+1}t^{b_j+1} \quad \text{for } j = 1, \dots, k, \\
    \lcm(s^{a_k+1}, t^{b_{k+1}+1}) &= s^{a_k+1}t^{b_{k+1}+1}.
\end{align*}
Since $s^{a_{j-1}+1}t^{b_j+1}$ divides $s^{a_{k+1}+1}t^{b_j+1}$ for all $j = 1, \dots, k$, the generators of the form $s^{a_{k+1}+1}t^{b_j+1}$ are redundant. 
Thus, we obtain
\[
\Ann(F_1,\dots,F_{k+1}) = (t^{b_1+1}, s^{a_1+1}t^{b_2+1}, \dots, s^{a_k+1}t^{b_{k+1}+1}, s^{a_{k+1}+1}).
\]
\end{proof}

To illustrate the algebraic structure underlying \Cref{thm:main}, we explicitly construct the graded minimal free resolution of the quotient ring $R/\Ann(F_1,\dots,F_r)$ where $F_i = S^{[a_i]}T^{[b_i]}$ for $i = 1,\dots,r$ with $a_1<\cdots<a_r$ and $b_1>\cdots>b_r$.
For notational convenience, let $a_0 := -1$ and $b_{r+1} := -1$. 
We denote the minimal generators of $\Ann(F_1, \dots, F_r)$ by $g_i := s^{a_{i-1}+1}t^{b_i+1}$ and $d_i := \deg(g_i)$ for $i = 1, \dots, r+1$.
Then, the graded minimal free resolution of $R/\Ann(W)$, $W := R\circ F_1 +\cdots+ R\circ F_r$ is given by the Hilbert-Burch complex of the form:
\[
0 \longrightarrow \bigoplus_{i=1}^r R[-a_i-b_i-2] \xrightarrow{\;\; M \;\;} \bigoplus_{i=1}^{r+1} R[-d_i] \longrightarrow R \longrightarrow R/\Ann(W) \longrightarrow 0,
\]
where the syzygy matrix $M$ is an $(r+1) \times r$ matrix defined as
\[
M = \begin{bmatrix}
    s^{a_1+1} & 0 & 0 & \dots & 0 \\
    -t^{b_1-b_2} & s^{a_2-a_1} & 0 & \dots & 0 \\
    0 & -t^{b_2-b_3} & s^{a_3-a_2} & \dots & 0 \\
    \vdots & \vdots & \ddots & \ddots & \vdots \\
    0 & 0 & \dots & -t^{b_{r-1}-b_r} & s^{a_r-a_{r-1}} \\
    0 & 0 & \dots & 0 & -t^{b_r+1}
\end{bmatrix}.
\]

Applying the canonical dualizing functor $\Hom_R(-, R[-2])$ to this resolution, we obtain the following dual resolution:
\[
0 \longrightarrow R[-2] \longrightarrow \bigoplus_{i=1}^{r+1} R[d_i-2] \xrightarrow{\;\; M^T \;\;} \bigoplus_{i=1}^r R[a_i+b_i] \longrightarrow \omega_{R/\Ann(W)} \longrightarrow 0.
\]
By apolarity duality, this dualized complex is precisely the graded minimal free resolution of $W = R \circ F_1 + \dots + R \circ F_r$.
In particular, each column vector of $M^T$ explicitly represents the algebraic relations among the generators of $W$ induced by the differential operators that annihilate the forms $F_1, \dots, F_r$.

Using \Cref{lem:monomialsimultaneousrank}, the following theorem provides the precise simultaneous and simultaneous border ranks of binary monomials.
\begin{theorem}\label{thm:monomialrank}
Suppose that $F_i = S^{a_i}T^{b_i}$ for $i = 1, \dots, r$ are binary monomials with $a_1 < \dots < a_r$ and $b_1 > \dots > b_r$.
For notational convenience, let $a_0 = b_{r+1} = -1$.
\begin{enumerate}
    \item\label{prop:monomialhonestrank} 
    If $a_1 > 0$ and $b_r > 0$, then the simultaneous rank of $F_1, \dots, F_r$ is exactly 
    \[
    \max\{b_1, a_r\} + 1.
    \]
    More precisely, it is equal to
    \[
    \begin{cases}
        b_1+1 &\quad\text{if}\quad a_r \leq b_1 \quad\text{and}\quad a_1 > 0, \\
        a_r+1 &\quad\text{if}\quad b_1 \leq a_r \quad\text{and}\quad b_r > 0.
    \end{cases}
    \]
    \item\label{prop:monomialborderrank} The simultaneous border rank of the monomials $F_1,\dots,F_r$ is 
    \[
    \min\{a_i+b_{i+1}+1 \mid 0\leq i\leq r\}+1.
    \]
\end{enumerate}
\end{theorem}
\begin{proof}\eqref{prop:monomialhonestrank} 
Let $d$ be any positive integer such that $1 \leq d \leq \max\{a_r, b_1\}$. 
We claim that the simultaneous rank is strictly larger than $d$.

Suppose $a_r \leq b_1$ and $a_1>0$. 
In this case, observe that 
\[
\Ann(F_1, \dots, F_r)_d = (t^{b_1+1},s^{a_1+1}t^{b_2+1},\dots,s^{a_r+1})_d = (s^{a_1+1}t^{b_2+1},\dots,s^{a_r+1})_d.
\]
Since $\gcd(s^{a_1+1}t^{b_2+1},\dots,s^{a_r+1}) = s^{a_1+1}$ and $a_1+1\geq 2$, by \Cref{thm:simultaneousrankalgorithm}, the existence of a non-square-free monomial factor in $\Ann(F_1,\dots,F_r)_d$ implies that the simultaneous rank must be strictly larger than $d$.

Suppose $a_r \geq b_1$ and $b_r>0$. 
By the symmetric argument, 
\[
\Ann(F_1,\dots,F_r)_d = (t^{b_1+1}, \dots, s^{a_{r-1}+1}t^{b_r+1}, s^{a_r+1})_d = (t^{b_1+1},\dots,s^{a_{r-1}+1}t^{b_r+1})_d.
\]
Since $\gcd(t^{b_1+1},\dots,s^{a_{r-1}+1}t^{b_r+1}) = t^{b_r+1}$ and $b_r+1\geq 2$ by \Cref{thm:simultaneousrankalgorithm}, the existence of a non-square-free monomial factor in $\Ann(F_1,\dots,F_r)_d$ implies that the simultaneous rank must be strictly larger than $d$.

Now, suppose that $d = \max\{b_1,a_r\}+1$.
Then the minimal monomial generators of $\Ann(F_1,\dots,F_r)$ of degree at most $d$ have greatest common divisor $1$ since they include $t^{b_1+1}$ and $s^{a_r+1}$.
Thus, by \Cref{thm:simultaneousrankalgorithm}, we conclude that the simultaneous rank is at most $\max\{b_1,a_r\}+1$.
By combining all arguments above, we obtain the desired equality.

\eqref{prop:monomialborderrank} 
The simultaneous border rank immediately follows from \Cref{cor:bordersimultaneousrank} and \Cref{lem:monomialsimultaneousrank} \eqref{lem:monomialAnn}.
\end{proof}
Under the conditions $a_1 > 0$ and $b_r > 0$, \Cref{thm:monomialrank} implies that the simultaneous border rank and the simultaneous rank coincide if and only if 
\[
\min\{a_i+b_{i+1}+1 \mid 0\leq i \leq r\} = b_1 = a_r,
\]
which is a highly restrictive case.

\section{Rational curves}\label{sec:rationalCurves}

Let 
\[
\mathcal{C} := \overline{\{[f_0(\alpha,\beta) : \dots : f_n(\alpha,\beta)] \mid [\alpha:\beta] \in \mathbb{P}_{\mathbb{K}}^1\}} \subset \mathbb{P}_{\mathbb{K}}^n
\]
be a non-degenerate projective curve parametrized by linearly independent binary forms $f_0, \dots, f_n \in \mathbb{K}[\alpha,\beta]_d$, where $\overline{\mathcal{X}}$ denotes the Zariski closure of a set $\mathcal{X}$. 
Note that every rational curve admits such a parameterization and the non-degeneracy of the curve implies $d \geq n$.

When $d = n$, the curve $\mathcal{C}$ is a rational normal curve of degree $d$. 
If $d > n$, it can be realized as the image of a degree-$d$ rational normal curve under a linear projection from a $(d-n-1)$-dimensional subspace. 
This can be deduced by combining standard results on projective rational maps (cf. \cite[Section 2]{MR3592466}).

More explicitly, regarding $V=\langle f_0,\dots,f_n\rangle$ as a subspace of $\mathbb{K}[\alpha,\beta]_d\cong R_d$, we take $L=V^\perp\subseteq\mathcal{R}_d$ with respect to the apolar pairing induced by the differentiation action.
Then $\mathcal{C}$ is projectively equivalent to the image of the degree-$d$ rational normal curve under the linear projection from the center $\mathbb{P}(L)$.

Now we reverse our viewpoint. 
Consider the image $\mathcal{C} = \overline{\pi_{\mathcal{L}}(\widetilde{\mathcal{C}}_d \setminus \mathcal{L})}$ of a rational normal curve $\widetilde{\mathcal{C}}_d \subset \mathbb{P}_{\mathbb{K}}^d$ under the linear projection $\pi_{\mathcal{L}} : \mathbb{P}_{\mathbb{K}}^d \dashrightarrow \mathbb{P}_{\mathbb{K}}^{d-r}$ where $\mathcal{L} \subset \mathbb{P}_{\mathbb{K}}^d$ is a linear space with $\dim \mathcal{L} = r-1$ for $1\leq r \leq d-1$. 
For this linear space $\mathcal{L}$, let $L$ be the $r$-dimensional vector subspace of $\mathbb{K}[S,T]_d$ corresponding to $\mathcal{L}$ (i.e., $\mathcal{L} = \mathbb{P}(L)$). 
Then, we have the rational map
\[
\pi_{\mathcal{L}} : \mathbb{P}(\mathbb{K}[S,T]_d) \dashrightarrow \mathbb{P}(\mathbb{K}[S,T]_d / L) \cong \mathbb{P}((L^\perp)^*) \cong \mathbb{P}_{\mathbb{K}}^{d-r}.
\]
Moreover, if we employ the pairing $\circ: \mathbb{K}[s,t]_d \times \mathbb{K}[S,T]_d \rightarrow \mathbb{K}$ introduced in \Cref{subsec:spaceDuality}, then $L^\perp = (\Ann(F_1,\dots,F_r))_d$ where $\{F_1, \dots, F_r\}$ is a basis for $L$.
\subsection{Apolar Defect} \label{subsec:ApolarDefect}
The structure of the annihilator $\Ann(L) = \Ann(F_1,\dots,F_r)$ is a special case of \Cref{thm:main} where the forms $F_1,\dots,F_r$ all have the same degree. 
The following proposition provides several equivalent conditions for this setting. 
\begin{proposition}
\label{prop:level_equivalences}
Let $I\subset R$ be a homogeneous Artinian ideal minimally generated by $r+1$ forms of degrees $d_1,\ldots,d_{r+1}$, and let
\begin{equation}\label{equ:Hilbert-BurchEquiDeg}
    0\longrightarrow \bigoplus_{j=1}^{r}R(-b_j)
\longrightarrow \bigoplus_{i=1}^{r+1}R(-d_i)
\longrightarrow R \longrightarrow R/I \longrightarrow 0
\end{equation}
be its minimal graded free resolution.
For a fixed integer $d\geq 0$, the following conditions are equivalent.
\begin{enumerate}[label=\textnormal{(\roman*)}]
    \item\label{item:level_equivalences1} There exists a $r$-dimensional subspace $L\subset \mathcal{R}_d$ such that $I=\Ann(L)$.
    \item\label{item:level_equivalences2} The algebra $A=R/I$ is level of socle degree $d$ and type $r$.
    \item\label{item:level_equivalences3} The inverse-system module
    \[
    I^{\perp}=\{F\in\mathcal{R}\mid I\circ F=0\}
    \]
    is minimally generated by $r$ linearly independent forms of degree $d$.
    \item\label{item:level_equivalences4} The last shifts satisfy
    \[
    b_1=\cdots=b_r=d+2.
    \]
    \item\label{item:level_equivalences5} One has
    \[
    \dim_{\mathbb{K}}I_d=d+1-r
    \qquad\text{and}\qquad
    I_{d+1}=R_{d+1}.
    \]
\end{enumerate}
\end{proposition}

\begin{proof}
The equivalence of \ref{item:level_equivalences1}, \ref{item:level_equivalences2}, and \ref{item:level_equivalences3} is the usual Macaulay inverse system correspondence: the minimal generators of $I^\perp$ are dual to the socle generators of $R/I$.
Since $R$ has two variables, a socle generator of degree $e$ corresponds to a summand $R(-e-2)$ in the last free module.  Hence \ref{item:level_equivalences2} is equivalent to \ref{item:level_equivalences4}.

If $A$ is level of socle degree $d$ and type $r$, then $A_{d+1}=0$ and $A_d=\operatorname{Soc}(A)$ has dimension $r$.  
This gives \ref{item:level_equivalences5}.
Conversely, assuming \ref{item:level_equivalences5}, one has $A_d\subseteq \operatorname{Soc}(A)$ and $\dim_{\mathbb{K}}A_d=r$.
The Hilbert--Burch resolution \eqref{equ:Hilbert-BurchEquiDeg} shows that the total socle dimension is $r$ because the last free module has rank $r$.
Therefore $\operatorname{Soc}(A)=A_d$, proving that $A$ is level of socle degree $d$ and type $r$.
\end{proof}

In particular, for $L\in\operatorname{Gr}(r,\mathcal{R}_d)$, the ideal $I=\Ann(L)$ is minimally generated by $r+1$ forms and has a resolution
\[
0\longrightarrow R(-d-2)^r
\longrightarrow \bigoplus_{i=1}^{r+1}R(-d_i)
\longrightarrow I\longrightarrow 0.
\]
Consequently,
\begin{equation}\label{equ:totalDegree}
    \sum_{i=1}^{r+1}d_i=r(d+2).
\end{equation}
\begin{definition}
The \emph{apolar defect} (or \emph{socle generation defect}) of $L = \langle F_1, \dots, F_r \rangle \subseteq \mathbb{K}[S,T]_d$, denoted by $\rho(L)$, is defined as the number of minimal generators of $\Ann(L)$ of degree $d+1$.
\end{definition}
The apolar defect $\rho(L)$ will play an important role in our subsequent geometric analysis. 
We investigate its properties and equivalent characterizations using known results on level algebras and ancestor ideals from \cite{MR2028071}.

For a subspace $L \subset \mathcal{R}_d$, we define its \emph{first prolongation} $L^{(1)} \subset \mathcal{R}_{d+1}$ as
\[
L^{(1)} = \{G \in \mathcal{R}_{d+1} \mid R_1 \circ G \subseteq L\} = \{G \in \mathcal{R}_{d+1} \mid s \circ G, t \circ G \in L\}.
\]
\begin{theorem}\label{thm:inverse_space_top_generators}
Let $L\in\operatorname{Gr}(r,\mathcal{R}_d)$ and $\Ann(L)\subset R = \mathbb{K}[s,t]$. 
Then
\[
\rho(L)=\dim_{\mathbb{K}}L^{(1)}.
\]
Equivalently, $\Ann(L)$ has a minimal generator of degree $d+1$ if and only if there exists a nonzero $G\in\mathcal{R}_{d+1}$ such that $s\circ G,t\circ G\in L$.
\end{theorem}
\begin{proof}
Let $I$ be the annihilator of $L$.
Since $I_{d+1}=R_{d+1}$, the space of minimal generators of $I$ in degree $d+1$ is
\[
\frac{I_{d+1}}{R_1I_d}=\frac{R_{d+1}}{R_1I_d}.
\]
Under the perfect apolar pairing in degree $d+1$, the dual of this quotient is $(R_1I_d)^\perp$.  
For $G\in\mathcal{R}_{d+1}$, one has
\[
G\in(R_1I_d)^\perp
\iff f\circ(\ell\circ G)=0
\quad\text{for every }f\in I_d\text{ and }\ell\in R_1.
\]
Since $I_d=L^\perp$, this is equivalent to $\ell\circ G\in L$ for every $\ell\in R_1$.  Thus $(R_1I_d)^\perp=L^{(1)}$.
\end{proof}
Let $J \subseteq R$ be the ancestor ideal of $I_d = \Ann(L)_d$ (see \cite{MR2028071} for its definition and properties).
Note that $J\cap R_{\geq d} = (I_d)$ and $J_{d+1}=R_1I_d$.
Then \Cref{thm:inverse_space_top_generators} implies that $L^{(1)}$ coincides with the apolar dual of $J_{d+1}$; that is, $L^{(1)} = J_{d+1}^\perp$.
As a result, we obtain
\[
\rho(L) = \dim_{\mathbb{K}} R_{d+1} - \dim_{\mathbb{K}}J_{d+1} = \dim_{\mathbb{K}}(J_{d+1}^\perp).
\]
In particular, the apolar defect is related to Iarrobino's invariant by
\[
\tau(L^\perp) = \dim_{\mathbb K}(R_1L^\perp)-\dim_{\mathbb K}L^\perp = r+1-\rho(L).
\]
By \cite[Lemma 2.2]{MR2028071}, this quantity equals the number of minimal generators of the ancestor ideal $J$.

Accordingly, let $g_1,\dots,g_m$ be the minimal homogeneous
generators of $\Ann(L)$ of degree at most $d$, where
$m=\tau(L^\perp)=r+1-\rho(L)$.
In what follows, we explore the relationship between the apolar defect $\rho(L)$, the greatest common divisor $\gcd(g_1, \dots, g_m)$, and the dimension of $L$.

\begin{corollary}\label{cor:inequalityApolarDefect}
The degree of $g=\gcd(g_1,\ldots,g_m)$ is bounded by the apolar defect of $L$; that is,
\[
\deg(g) \leq \rho(L) \leq \dim_{\mathbb{K}}(L).
\]
\end{corollary}

\begin{proof}
Let $d'=\deg(g)$ and $J=(g_1,\ldots,g_m)$.
Then $J_d\subset gR_{d-d'}$. 
Since $J$ is the ancestor ideal of $I$, by \Cref{prop:level_equivalences}\ref{item:level_equivalences5}, $\dim_{\mathbb{K}}J_d = \dim_{\mathbb{K}} \Ann(L)_d = d+1-r$ and $\dim_{\mathbb{K}}J_{d+1}=d+2-\rho(L)$.
This yields 
\[
(d+2)-\rho(L) \leq d-d'+2 \implies d' \leq \rho(L).
\]
Moreover, $\Ann(L)_d \not= 0$ since $r \leq d - 1$, so $m\geq 1$.
Thus,
\[
\rho(L) = r+1-m \leq r = \dim_{\mathbb{K}}(L).
\]
\end{proof}

As we will show in \Cref{subsec:Parametrization}, the apolar defect $\rho(L)$ not only dictates the parametrization of the curve and the dimension of the containing scroll, but also bounds the length of the intersection scheme $\mathcal{L} \cap \widetilde{\mathcal{C}}_d$.

When $L$ is spanned by monomials (i.e., $\mathcal{L}$ is a coordinate linear subspace), its apolar defect can be explicitly computed by counting consecutive exponents.

\begin{corollary}\label{cor:d+1gensForMonomials}
Let $L \subset \mathbb{K}[S,T]_d$ be the subspace spanned by monomials $\{S^{d-c_i}T^{c_i} \mid 1 \le i \le r\}$ with $0 \le c_1 < \dots < c_r \le d$.
Then
\[
\rho(L) = |\{ 1 \le i \le r-1 \mid c_{i+1} - c_i = 1 \}| + |\{c_1, c_r\} \cap \{0, d\}|.
\]
\end{corollary}

\begin{proof}
This follows immediately from \Cref{lem:monomialsimultaneousrank}\eqref{lem:monomialAnn}.
\end{proof}

\Cref{cor:d+1gensForMonomials} clarifies why the inequality $\deg(g) \le \rho(L)$ holds and how a gap between the two terms arises. For linear spaces $L$ spanned by monomials, $\deg(g)$ is obtained by summing the lengths of the maximal consecutive exponent blocks containing the boundary points $0$ (if $c_1 = 0$) and $d$ (if $c_r = d$). 
In contrast, $\rho(L)$ accounts for all adjacent pairs of consecutive integers across the entire sequence $(c_1, \dots, c_r)$, including internal ones, along with the endpoints in $\{0, d\}$.

\subsection{Parametrization}\label{subsec:Parametrization}
We obtain parametrizations of rational curves up to projective equivalence.
\begin{theorem}\label{thm:parametrizedCurve}
Let $L = \langle F_1, \dots, F_r \rangle \subseteq \mathcal{R}_d$ have dimension $r$, and $m := r + 1 - \rho(L)$. 
Let $\mathcal{L} = \mathbb{P}(L)$ be the linear subspace in $\mathbb{P}_{\mathbb{K}}^d = \mathbb{P}(\mathcal{R}_d)$.
Let $g_1, \dots, g_m$ be the minimal generators of $\Ann(L)$ of degree at most $d$ with $d_i := \deg(g_i)$, and let $g = \gcd(g_1, \dots, g_m)$ with $g_i' := g_i/g$. 
Then the projected curve $\mathcal{C} = \overline{\pi_{\mathcal{L}}(\widetilde{\mathcal{C}}_d \setminus \mathcal{L})}$ is projectively equivalent to the curve parametrized by \eqref{equ:parametrizationCurve}.
\end{theorem}
\begin{proof}
Given the $r$-dimensional linear space $L \subset \mathbb{K}[S,T]_d$, we have 
\[
L^\perp = \left(\operatorname{Ann}(L)\right)_d = \bigoplus_{i=1}^m g_i \mathbb{K}[s,t]_{d-d_i}.
\]
Note that this direct sum decomposition can be verified by a straightforward dimension-counting argument on both sides.

Since $\dim L^\perp = d+1-r$, the following $d-r+1$ polynomials form a basis for $L^\perp \subseteq \mathbb{K}[s,t]_d$:
\[
\begin{array}{cccc}
v_r := g_m s^{d-d_m}, & v_{r+1} := g_m s^{d-d_m-1}t, & \dots, & v_{r+d-d_m} := g_m t^{d-d_m}, \\[0.3em]
\vdots & \vdots & & \vdots \\[0.3em]
v_{d_1} := g_1 s^{d-d_1}, & \dots, & v_{d-1} := g_1 s t^{d-d_1-1}, & v_d := g_1 t^{d-d_1}
\end{array}
\]

We consider $\mathbb{K}[S,T]_d$ as the dual vector space to $\mathbb{K}[s,t]_d$ via the non-degenerate apolar product $\circ$. We choose $\{w_r, \dots, w_d\} \subset \mathbb{K}[S,T]_d$ as a dual system to $\{v_r, \dots, v_d\}$, satisfying $v_i \circ w_j = \delta_{ij}$ for $i,j \in \{r, \dots, d\}$. 
Note that $L \cap \langle w_r, \dots, w_d \rangle = \{0\}$ because $v \circ F = 0$ for all $v \in L^\perp$ and $F \in L$. 
Thus, we can extend this to a basis $\{v_0, \dots, v_d\}$ of $\mathbb{K}[s,t]_d$ which is dual to a basis $\{F_1, \dots, F_r, w_r, \dots, w_d\}$ of $\mathbb{K}[S,T]_d$.

For any point $q = [F_q] \in \mathbb{P}_{\mathbb{K}}^d$, we have the expansion:
\[
F_q = \sum_{j=1}^r \alpha_j F_j + \sum_{i=r}^d (v_i \circ F_q) w_i.
\]
The projection $\pi_{\mathbb{P}(L)}: \mathbb{P}_{\mathbb{K}}^d \setminus \mathbb{P}(L) \to \mathbb{P}_{\mathbb{K}}^{d-r}$ maps $q$ to the coefficients of the $w_i$ components:
\[
\pi_{\mathbb{P}(L)}(q) = [v_r \circ F_q : v_{r+1} \circ F_q : \dots : v_d \circ F_q].
\]
Let $q = [\alpha^d : \alpha^{d-1}\beta : \dots : \beta^d] \in \widetilde{\mathcal{C}_d}$ be a point on the rational normal curve, and choose the representative $F_q = (\alpha S + \beta T)^d$. 
By \Cref{lem:evaluation}, we have
\begin{equation*}
    (g_i s^{d-d_i-k} t^k) \circ (\alpha S + \beta T)^d = d! g_i(\alpha, \beta) \alpha^{d-d_i-k} \beta^k.
\end{equation*}
By dividing each coordinate by the greatest common divisor $g$, we obtain the desired parametrization.
\end{proof}
The parametrization of $\mathcal{C}$ in \eqref{equ:parametrizationCurve} also gives a geometric interpretation of the common factor $g$: it records, with multiplicity, the points responsible for the inner projections involved in $\pi_{\mathcal L}$.
Taking into account the generic degree of the parametrization, we also obtain a formula for the degree of the projected curve $\mathcal{C}$ as follows.
\begin{corollary}\label{cor:detectInnerProjections}
With the same notation as in \Cref{thm:parametrizedCurve}, set
$d'=\deg(g)$.
Let $\varphi:\mathbb{P}^1_{\mathbb{K}}\longrightarrow\mathcal{C}$ denote the morphism defined by the parametrization \eqref{equ:parametrizationCurve}.
Let $\delta$ be its generic degree onto its image; that is, a general point of $\mathcal{C}$ has $\delta$ inverse images under $\varphi$, counted with multiplicity.
Then
\[
\deg(\mathcal{C})=\frac{d-d'}{\delta}
\quad\text{and}\quad
\operatorname{length}
\bigl(\mathcal L\cap\widetilde{\mathcal{C}}_d\bigr)=d'.
\]
\end{corollary}

\begin{proof}
For $1\leq i\leq m$ and $0\leq j\leq d-d_i$, let $h_{i,j} = g_i'(\alpha,\beta)\alpha^{d-d_i-j}\beta^j$.
By \Cref{thm:parametrizedCurve}, the parametrization of
$\mathcal{C}$ is given by the forms $h_{i,j}$, all of which have degree $d-d'$.
Note that, since $\gcd(g_1',\dots,g_m')=1$, the forms $h_{i,j}$ have no common zero on $\mathbb{P}^1_{\mathbb{K}}$.

We first claim that $\delta\deg(\mathcal{C})=d-d'$.
Let $H\subset\mathbb{P}^{d-r}_{\mathbb{K}}$ be a general hyperplane.
Since $\varphi$ has generic degree $\delta$, the inverse image of $H\cap\mathcal{C}$ under $\varphi$ has length $\delta\deg(\mathcal{C})$. 
On the other hand, this inverse image is the vanishing scheme of a nonzero linear combination of the forms $h_{i,j}$ and hence has length $d-d'$. 
Therefore, $\delta\deg(\mathcal{C})=d-d'$.

Now, we compute the length of $\mathcal L\cap\widetilde{\mathcal{C}}_d$.
Via the apolar pairing, the elements $g_i(s,t)s^{d-d_i-j}t^j\in L^\perp$ where $1\leq i\leq m$ and $0\leq j\leq d-d_i$ represent the linear equations defining $\mathcal L$.
Under the standard parametrization $q=[\alpha^d:\alpha^{d-1}\beta:\cdots:\beta^d]$ of $\widetilde{\mathcal{C}}_d$, the pullbacks of these linear equations are the degree-$d$ forms $g_i(\alpha,\beta)\alpha^{d-d_i-j}\beta^j=g h_{i,j}$.
Thus, the intersection scheme $\mathcal L\cap\widetilde{\mathcal{C}}_d$ corresponds to the common zero scheme of the forms $g h_{i,j}$.

Since the forms $h_{i,j}$ have no common zero, the common zeros of the forms $g h_{i,j}$ are precisely the zeros of $g$ (with the same multiplicities). 
Thus we have $\operatorname{length} \bigl(\mathcal L\cap\widetilde{\mathcal{C}}_d\bigr)=d'$.
\end{proof}

\begin{corollary}\label{cor:RNC}
Assume, in addition, that the morphism $\varphi:\mathbb{P}_{\mathbb{K}}^1 \rightarrow \mathcal{C}$ is birational onto its image.
Then the following are equivalent:
\begin{enumerate}[label=\textnormal{(\roman*)}]
    \item\label{item:inequalityApolarDefect1} $\deg(g) = \dim_{\mathbb{K}}(L)$.
    \item\label{item:inequalityApolarDefect2} $\rho(L) = \dim_{\mathbb{K}}(L)$.
    \item\label{item:inequalityApolarDefect3} $\mathcal{C} = \overline{\pi_{\mathcal{L}}(\widetilde{\mathcal{C}_d}\setminus \mathcal{L})}\subset \mathbb{P}_{\mathbb{K}}^{d-r}$ is a rational normal curve of degree $d-r$.
\end{enumerate}
\end{corollary}
\begin{proof}
For the equivalences, the implication \ref{item:inequalityApolarDefect1} $\implies$ \ref{item:inequalityApolarDefect2} follows immediately from the above inequalities, and \ref{item:inequalityApolarDefect2} $\implies$ \ref{item:inequalityApolarDefect3} holds by \Cref{thm:parametrizedCurve}.
To prove \ref{item:inequalityApolarDefect3} $\implies$ \ref{item:inequalityApolarDefect1}, suppose that $\mathcal{C}$ is a rational normal curve of degree $d-r$.
Since $\varphi$ is birational onto its image, \Cref{cor:detectInnerProjections} gives $\deg(g)=d-\deg(\mathcal{C})=r=\dim_{\mathbb{K}}(L)$.
\end{proof}
For one-point outer projections with $d\geq3$, \cite[Theorem 7.3]{MR2274517} gives a containing rational normal surface scroll.
In the range of \Cref{thm:minimalSectionScroll}, our parametrization constructs a containing scroll of dimension $r+1-\rho(L)$, and the following argument proves its minimality under the section condition.

Let $p:\mathbb P(\bigoplus_{i=1}^n \mathcal{O}_{\mathbb{P}^1_{\mathbb{K}}}(a_i)) \to\mathbb{P}^1_{\mathbb{K}}$ be the natural projection.
For each $q\in\mathbb P^1_{\mathbb K}$, the fiber of $p$ parametrizes one-dimensional quotients of $\bigoplus_{i=1}^n\mathcal{O}_{\mathbb{P}^1_{\mathbb{K}}}(a_i)|_q$.
Let $j$ denote the morphism from the projective bundle $\mathbb P(\bigoplus_{i=1}^n \mathcal{O}_{\mathbb{P}^1_{\mathbb{K}}}(a_i))$ to $\mathcal{S}=S(a_1,\dots,a_n)$ defined by the global sections of its tautological line bundle.
We say that \emph{$\mathcal{C}$ is contained as a section of its ruling}
if there exists a morphism
$\sigma:\mathbb P^1_{\mathbb K}\to
\mathbb P(\bigoplus_{i=1}^n\mathcal O_{\mathbb P^1_{\mathbb K}}(a_i))$
such that $p\circ\sigma=\operatorname{id}_{\mathbb P^1_{\mathbb K}}$
and the induced morphism
$\psi=j\circ\sigma:\mathbb P^1_{\mathbb K}\to\mathcal S$
has image $\mathcal C$.
With this convention, the following proof shows that, in the range of \Cref{thm:minimalSectionScroll},
$r+1-\rho(L)$ is the minimum dimension of a rational normal
scroll containing $\mathcal{C}$ as a section of its ruling.

\begin{proof}[Proof of \Cref{thm:minimalSectionScroll}]
Let $L^\perp=\Ann(L)_d=gW$ where $W\subseteq R_{d-d'}$ is the linear system defining $\varphi$.
Then $\dim_{\mathbb{K}}W = d-r+1$.
Since $\Ann(L)_{d+1}=R_{d+1}$, 
we have 
\[
\rho(L)=d+2-\dim_{\mathbb{K}}(R_1L^\perp)=d+2-\dim_{\mathbb{K}}(R_1W).
\]
Suppose first that $\rho(L)=r+1-m$.
Then $\Ann(L)$ has exactly $m$ minimal generators of degree at most $d$.
Let these generators be $g_1,\dots,g_m$, with $d_i=\deg(g_i)$ and $g_i'=g_i/g$, as in \Cref{thm:parametrizedCurve}.
Then $W=\bigoplus_{i=1}^m g_i'R_{d-d_i}$.
By taking dimensions, we obtain 
\[
\sum_{i=1}^m(d-d_i)=d-r+1-m>0.
\]
Therefore $d-d_i>0$ for some $i$.
For this $i$, since $\frac{g_i'\alpha^{d-d_i-1}\beta}{g_i'\alpha^{d-d_i}}=\frac{\beta}{\alpha}$ on a nonempty open subset, $\varphi$ is birational onto its image.
That is, $\delta = 1$.

By \Cref{thm:parametrizedCurve}, after a projective change of coordinates, $\mathcal{C}$ is parametrized by \[[\dots:g_i'(\alpha,\beta)\alpha^{d-d_i-j}\beta^j:\dots]\] where $1\leq i\leq m$ and $0\leq j\leq d-d_i$.
Consider the rational normal scroll $\mathcal S=S(d-d_1,\dots,d-d_m)\subseteq\mathbb{P}_{\mathbb{K}}^{d-r}$ with the following parametrization:
\begin{equation}
\{[\dots:\gamma_i\eta^{d-d_i-j}\zeta^j:\dots]_{1\le i \le m, 0\le j \le d-d_i}\mid [\eta:\zeta]\in \mathbb{P}_{\mathbb{K}}^1, [\gamma_1:\cdots:\gamma_m]\in\mathbb{P}_{\mathbb{K}}^{m-1}\}.
\end{equation}
Since $d-d_i\geq 0$ for every $i = 1,\dots,m$ and $\sum_{i=1}^m (d-d_i) >0$, this scroll has dimension $m$.
Observe that the parametrization of $\mathcal{C}$ is recovered by setting $[\eta:\zeta]=[\alpha:\beta]$ and, in local ruling coordinates, $\gamma_i=g_i'(\alpha,\beta)$.
Since the forms $g_i'$ have no common zero, this defines a section of the ruling family over $\mathbb{P}_{\mathbb{K}}^1$ whose image is $\mathcal{C}$.
Thus $\mathcal{C}$ is contained in $\mathcal S$ as a section of its ruling (understood on the standard projective-bundle model when $\mathcal S$ is a cone).

To prove minimality of $m$, suppose that $\varphi$ is birational and that $\mathcal{C}$ is contained in an $n$-fold rational normal scroll $\mathcal S=S(a_1,\dots,a_n) \subseteq \mathbb{P}_{\mathbb{K}}^{d-r}$ as a section of its ruling.
Then $\sum_{i=1}^n(a_i+1)=d-r+1$.
Let $\psi:\mathbb P^1_{\mathbb K}\to\mathcal{C}$ be the morphism
induced by the section.
We claim that $\psi$ is birational onto $\mathcal{C}$.
Indeed, since $\mathcal{C}$ is nondegenerate, it is not contained in the vertex of $\mathcal{S}$.
Away from the vertex, $j$ is an isomorphism, and $p\circ j^{-1}$
gives an inverse to $\psi$ on a nonempty open subset of $\mathcal{C}$.

Since $\varphi$ and $\psi$ are birational morphisms from the nonsingular projective curve $\mathbb{P}^1_{\mathbb{K}}$, they are finite and hence are normalization morphisms of $\mathcal{C}$ \cite[Chapter II, Proposition~6.8 and Exercise~3.8]{hartshorne}.
The uniqueness of normalization therefore yields a unique $\sigma\in\operatorname{Aut}_{\mathbb{K}}(\mathbb{P}^1_{\mathbb{K}})$ such that $\psi\circ\sigma=\varphi$.

After identifying the base of the ruling with the domain of $\varphi$ via $\sigma$, the section parametrization agrees with $\varphi$.
The section corresponds to a surjection
\[
\bigoplus_{i=1}^n
\mathcal O_{\mathbb{P}_{\mathbb K}^1}(a_i)
\twoheadrightarrow
\varphi^*\mathcal O_{\mathcal C}(1)
\cong
\mathcal O_{\mathbb{P}_{\mathbb K}^1}(d-d')
\]
\cite[Chapter II, Proposition~7.12]{hartshorne}.
Its components are therefore forms $f_i\in R_{d-d'-a_i}$ with no common zero.
In homogeneous coordinates adapted to the scroll, the coordinate forms of $\varphi$ are $f_i\alpha^{a_i-j}\beta^j$, where $1\le i\le n$ and $0\le j\le a_i$.
Since a projective change of ambient coordinates only changes the basis of $W$, we obtain $W=\sum_{i=1}^n f_iR_{a_i}$.

Therefore, $R_1W=\sum_{i=1}^n f_iR_{a_i+1}$, and hence $\dim_{\mathbb{K}}(R_1W)\leq\sum_{i=1}^n(a_i+2)=d-r+n+1$.
It follows that 
\begin{equation}\label{equ:defectBound}
\rho(L)=d+2-\dim_{\mathbb{K}}(R_1W)\geq r+1-n.
\end{equation}
In particular, if $\rho(L)=r+1-m$, then $m\leq n$.
Together with the construction above, this proves that $m$ is the smallest possible dimension.

Conversely, suppose that $\varphi$ is birational and that $m$ is the smallest positive integer such that $\mathcal{C}$ is contained in an $m$-fold rational normal scroll as a section of its ruling.
Applying \eqref{equ:defectBound} with $n=m$, we obtain $\rho(L)\geq r+1-m$.
Suppose that this inequality is strict, and set $n=r+1-\rho(L)$.
Then $\Ann(L)$ has exactly $n$ minimal generators of degree at most $d$ with $1\leq n<m\leq d-r$.
Applying the same construction to these generators gives an $n$-fold rational normal scroll $\mathcal{S}$ containing $\mathcal{C}$ as a section of its ruling.
This contradicts the minimality of $m$.
Therefore, $\rho(L)=r+1-m$, as desired.
\end{proof}
Assuming that $\varphi$ is birational, $m = 1$ if and only if $\mathcal{C}\subset \mathbb{P}_{\mathbb{K}}^{d-r}$ is a rational normal curve of degree $d-r$ by \Cref{cor:RNC}.
Thus, if $d\geq 2r+1$, then \Cref{thm:minimalSectionScroll} covers the full range $1\leq m\leq r+1 = \min\{r+1,d-r\}$.

Together with this case, \Cref{thm:minimalSectionScroll} identifies $m=r+1-\rho(L)$ as the minimum dimension of a rational normal scroll (possibly a cone) containing $\mathcal{C}$ as a section of its ruling.
However, if the section condition is dropped, $\mathcal{C}$ may be contained in a rational normal scroll of dimension smaller than $m$.

\begin{remark}
Assume that $\varphi$ is birational onto its image and $d-r\geq3$.
By \Cref{cor:RNC}, \Cref{thm:minimalSectionScroll}, and Green's $\mathcal K_{p,1}$ theorem \cite{MR739785} (see also \cite[Theorem~1.1]{MR4816776}), the minimum dimension of a rational normal scroll containing $\mathcal{C}$ without imposing the section condition can be determined when $\rho(L)\geq r-2$:
\begin{enumerate}
    \item If $\rho(L)=r$, the minimum is $1$.
    \item If $\rho(L)=r-1$, the minimum is $2$.
    \item If $\rho(L)=r-2$ and $d-r\neq5$, the minimum is $2$ if $\mathcal{K}_{d-r-2,1}(\mathcal{C})\neq0$ and $3$ otherwise.
\end{enumerate}
In the third case with $d-r=5$, the same criterion holds unless $\mathcal{C}$ lies on a Veronese surface, in which case the minimum is $3$ (see \cite[Theorem~5.2]{MR4816776}).
\end{remark}
\begin{example}
For an integer $d \ge 6$, let $\mathcal{L} = \langle e_1, e_3 \rangle \subset \mathbb{P}^d$ be the line corresponding to $L = \langle S^{d-1}T, S^{d-3}T^3 \rangle \subset \mathbb{K}[S,T]_d$. 
The defining ideal of the projected curve $\mathcal{C} := \overline{\pi_{\mathcal{L}}(\widetilde{\mathcal{C}}\setminus \mathcal{L})} \subset \mathbb{P}^{d-2}$ is obtained by eliminating $x_1, x_3$ from $I(\widetilde{\mathcal{C}})$. 
Since $\operatorname{Ann}(L) = (t^4, s^{d-2}t^2, s^d)$, $\mathcal{C}$ is contained in the $3$-fold scroll $S(0,0,d-4)$ by \Cref{thm:minimalSectionScroll}. 
In fact, $\mathcal{C}$ lies as a divisor on the surface scroll $\mathcal{S}=S(\lceil(d-6)/2\rceil,\lfloor d/2\rfloor)$.
The corresponding curve on the standard projective-bundle model of $\mathcal{S}$ has bidegree $(2,5)$ for odd $d$ and $(2,6)$ for even $d$ with respect to the minimal section and fiber classes.
This embedding can be verified via the $1$-dimensional toric polytope of $\mathcal{C}$ or Green's $\mathcal{K}_{p,1}$ theorem \cite{MR739785}.
Similar phenomena occur for $\mathcal{L} = \langle e_1, e_{d-1} \rangle$ and $\mathcal{L} = \langle e_{d-3}, e_{d-1} \rangle$.
\end{example}

\subsection{Determinantal loci on the Grassmannian}\label{subsec:Antiderivatives}
As we can regard the $r$-dimensional subspace $L\subset \mathcal{R}_d$ as a point in the Grassmannian $\Gr(r,\mathcal{R}_d)$, we consider the locus of such linear spaces in the Grassmannian.
For an integer $\tau\ge0$, let
\[
\mathfrak D_{\ge\tau}
:=\{L\in\Gr(r,\cR_d)\mid\rho(L)\ge\tau\}.
\]
Then $\mathfrak{D}_\tau=\mathfrak{D}_{\ge\tau}\backslash \mathfrak{D}_{\ge\tau+1}$ is the locus of $L$ where $\Ann(L)$ has exactly $\tau$ degree $d+1$ minimal generators. 
As mentioned in the introduction, we will show that $\mathfrak{D}_{\geq \tau}$ is a Zariski closed set.

\begin{proposition}\label{prop:apolarDefectLocus}
For $\tau\geq 0$, 
\[
\mathfrak D_{\ge\tau} = \{L\in\Gr(r,\cR_d)\mid\rho(L)\ge\tau\}
\]
is a Zariski closed set in the Grassmannian. 
And it is defined by minors of a matrix.
\end{proposition}
\begin{proof}
Let $D_s^{(m)}$ and $D_t^{(m)}$ be differentiation matrices using normalized bases from $\cR_m$ to $\cR_{m-1}$. 
In fact, they are $m\times (m+1)$ matrices given by $\begin{bmatrix} I_{m} \mid \mathcal{O} \end{bmatrix}$ and $\begin{bmatrix} \mathcal{O} \mid I_{m} \end{bmatrix}$, respectively, where $\mathcal{O} \in \mathbb{K}^{m}$ is the zero column vector.
Using the Stiefel coordinates for the Grassmannian, represent $L$ by a full-rank matrix $A\in\operatorname{Mat}_{r\times(d+1)}(\mathbb{K})$ whose rows form a basis of $L$. 
Then $G\in L^{(1)}$ if and only if there exist $a,b\in \mathbb{K}^r$ such that
\[
D_s^{(d+1)}z=A^Ta, \quad \text{and} \quad D_t^{(d+1)}z=A^Tb,
\]
where $z$ is the coefficient vector of $G$. 
Define
\[
\mathcal N(A)=
\begin{pmatrix}
D_s^{(d+1)}&-A^T&0\\
D_t^{(d+1)}&0&-A^T
\end{pmatrix}.
\]
Then $G\in L^{(1)}$ if and only if \[\mathcal{N}(A)\begin{pmatrix} z\\a\\b\end{pmatrix}=\begin{pmatrix}D_s^{(d+1)}z-A^Ta\\D_t^{(d+1)}z-A^Tb\end{pmatrix}=\begin{pmatrix}0\\0\end{pmatrix}\] has a solution. 
So the projection onto $z$ gives $\ker\mathcal N(A)\cong L^{(1)}$, hence
\[
\dim_{\mathbb{K}} L^{(1)}=d+2+2r-\rank\mathcal N(A).
\]
Therefore,
\[
\dim_{\mathbb{K}} L^{(1)}\ge\tau \iff \rank\mathcal N(A)\le d+2+2r-\tau.
\]
The corresponding Stiefel equations are the $(d+3+2r-\tau)$-minors.
\end{proof}
We will introduce several ways to describe this locus.
\begin{remark}
There are several ways to understand this locus.
\begin{enumerate}
\item Choose a full-rank matrix $Q\in\operatorname{Mat}_{(d+1-r)\times(d+1)}(\mathbb{K})$ satisfying $QA^T=0$.
Then $Q$ represents $\cR_d\to\cR_d/L$ and, dually, the point $L^\perp\in\Gr(d+1-r,R_d)$. 
Define
\[
M(Q)=
\begin{pmatrix}
QD_s^{(d+1)}\\
QD_t^{(d+1)}
\end{pmatrix}.
\]
Then $L^{(1)}=\ker M(Q)$, and therefore $\dim_K L^{(1)}=d+2-\rank M(Q)$.
Thus
\[
\dim_KL^{(1)}\ge\tau \iff \rank M(Q)\le d+2-\tau,
\]
and the equations are the $(d+3-\tau)$-minors. 
On the standard chart $A=(I_r\mid X)$, one may take $Q=(-X^T\mid I_{d+1-r})$.
Moreover,
\[
\rank \mathcal{N}(A) = 2r + \rank M(Q).
\]
Actually, using appropriate column and row operations, we can change $\mathcal{N}(A)$ to
\[\begin{pmatrix} I_{2r}&0\\0&M(Q)\end{pmatrix}\]

\item Define the gradient map $\nabla: \mathcal{R} \to \mathcal{R}\oplus \mathcal{R}$ by $\nabla(G) = (s \circ G, t \circ G)$, and the curl map $\curl: \mathcal{R}\oplus\mathcal{R} \to \mathcal{R}$ by $\curl(P,Q) = t \circ P - s \circ Q$.
Then the homogeneous gradient--curl sequence
\[
0\longrightarrow\cR_{d+1}
\xrightarrow{\nabla}\cR_d^{\oplus2}
\xrightarrow{\curl}\cR_{d-1}
\longrightarrow0
\]
is exact.
Indeed, if $t\circ P=s\circ Q$, then $G=\frac{SP+TQ}{d+1}$ satisfies $s\circ G=P$ and $t\circ G=Q$ by Euler's identity $S(s\circ F)+T(t\circ F)=dF$.

The restriction of curl to $L\oplus L$ is represented by
\[
\Psi(A)=
\begin{pmatrix}
D_t^{(d)}A^T&-D_s^{(d)}A^T
\end{pmatrix}.
\]
Hence $L^{(1)} \cong \ker \Psi(A)$ which implies that $\dim_{\mathbb{K}} L^{(1)} = 2r - \rank \Psi(A)$.
Therefore
\[
\dim_{\mathbb{K}}L^{(1)} \geq \tau \iff \rank\Psi(A) \leq 2r - \tau,
\]
and the equations are the $(2r- \tau + 1)$-minors.

\item Combining the formulas gives
\[
\begin{aligned}
\rho(L)=\dim_\mathbb{K}L^{(1)}
&=d+2+2r-\rank\mathcal N(A)\\
&=d+2-\rank M(Q)\\
&=2r-\rank\Psi(A)\\
\end{aligned}
\]
\end{enumerate}
\end{remark}

\section{Rational curves on surface scrolls}\label{sec:curvesinSurfacescroll}
Once the parametrization $\mathcal{C} = \{[f_0(\alpha,\beta) : \dots : f_{d-r}(\alpha,\beta)] \mid [\alpha:\beta]\in \mathbb{P}_{\mathbb{K}}^1\} \subset \mathbb{P}^{d-r}$ is established, its defining equations can be obtained algebraically by eliminating the parameters $\alpha$ and $\beta$ from the ideal $(z_0 - f_0, \dots, z_{d-r} - f_{d-r})$.
That is, 
\[
I(\mathcal{C}) = (z_0 - f_0, \dots, z_{d-r} - f_{d-r}) \cap \mathbb{K}[z_0,\dots,z_{d-r}].
\]

However, exploiting the geometric properties of $\mathcal{C}$ often allows for a more efficient computation of its defining ideal. 
We focus on the case where the curve lies as a divisor in a $2$-fold rational normal scroll.

\subsection{Smooth rational curve contained in a smooth surface scroll}\label{subsec:smoothcurveinSurface}
Let $\mathcal{C}=\overline{\pi_{\mathcal L}(\widetilde{\mathcal{C}}_d\setminus \mathcal{L})}\subset\mathbb{P}_{\mathbb{K}}^{d-r}$ be the projected rational curve and $d'=\operatorname{length}(\mathcal L\cap\widetilde{\mathcal{C}}_d)$ and $r=\dim_{\mathbb{K}}L$.
Let $\varphi:\mathbb{P}_{\mathbb{K}}^1\to\mathcal{C}$ be the morphism defined by the parametrization in \eqref{equ:parametrizationCurve} and let $\delta=\deg(\varphi)$ denote its generic degree onto $\mathcal{C}$.
Recall that, by \Cref{cor:detectInnerProjections}, $\deg(\mathcal{C})=\frac{d-d'}{\delta}$.
The embedding of $S(a,b)$ in $\mathbb{P}_{\mathbb{K}}^{d-r}$ with $1 \le a \le b$ is given by the complete linear system $|C_0 + bF|$, where $C_0$ denotes the directrix section with $C_0^2 = -e = a-b$, and $F$ represents the fiber class.
Since the scroll $S(a,b)$ contains the nondegenerate curve $\mathcal{C}$, we have $a+b+1 = d-r$.
\begin{theorem}\label{thm:surface-stable}
Assume that the curve $\mathcal{C}$ is not a rational normal curve, but lies on a smooth rational normal surface scroll $S(a,b) \subset \mathbb{P}_{\mathbb{K}}^{d-r}$ with $1 \le a \le b$.
If $\mathcal{C}$ is a section of the ruling of $S(a,b)$, then $\mathcal{C}$ is smooth. 
Conversely, if $\mathcal{C}$ is smooth and $d\geq \frac{2}{3}(\deg(\mathcal{C})+1) +r$, then $\mathcal{C}$ is a section of the ruling of $S(a,b)$.
\end{theorem}
\begin{proof}
Let $\mathcal{C} \subset \mathcal{S} = S(a,b)$ be a curve of bidegree $(m, n)$ in $\mathrm{Pic}(\mathcal{S}) \cong \mathbb{Z}C_0 \oplus \mathbb{Z}F$. 
By applying the adjunction formula for effective divisors on smooth surfaces (cf. \cite[Chapter V, Ex.~1.3]{hartshorne}), we have
\begin{equation}\label{equ:adjunctionFormula}
    2p_a(\mathcal{C}) - 2 = \mathcal{C}.(\mathcal{C} + K_{\mathcal{S}}),
\end{equation}
where $K_{\mathcal{S}}$ is the canonical divisor of $\mathcal{S}$. 
Substituting $\mathcal{C} \equiv mC_0 + nF$ and $K_{\mathcal{S}} \equiv -2C_0 - (e+2)F$ into \eqref{equ:adjunctionFormula} and using $C_0^2 = -e$, $C_0.F = 1$, and $F^2 = 0$, we obtain
\begin{equation}\label{equ:genus}
  p_a(\mathcal{C}) = (m-1)(n-1) - e\binom{m}{2}.  
\end{equation}
In particular, if $m = 1$ (i.e., $\mathcal{C}$ is a section of the ruling of $\mathcal{S}$), we obtain $p_a(\mathcal{C}) = 0$.

Conversely, suppose $\mathcal{C}$ is smooth and $2\deg(\mathcal{C}) \leq 3(d-r)-2$.
By \eqref{equ:genus}, we have $(m-1)(2n-em-2) = 0$.
Assume for contradiction that $\mathcal{C}$ is not a section meaning that $m \neq 1$. 
Then we obtain
\begin{equation}\label{equ:smoothcondition1}
    2n = (b-a)m+2.
\end{equation}
Since the hyperplane section class of $S(a,b)$ is $H = C_0+bF$, its degree is given by 
\begin{equation}\label{equ:smoothcondition2}
\deg(\mathcal{C}) = H.\mathcal{C} = (C_0+bF).(mC_0+nF) = am+n. 
\end{equation}
Combining \eqref{equ:smoothcondition1} and \eqref{equ:smoothcondition2} yields 
\[
m = \frac{2\deg(\mathcal{C})-2}{a+b} = 2+\frac{2(\deg(\mathcal{C})-d+r)}{d-r-1}.
\]
Since $\mathcal{C}$ is not a rational normal curve, $\deg(\mathcal{C}) > d-r$.
Moreover, the condition $d\geq \frac{2}{3}(\deg(\mathcal{C})+1) +r$ implies $2(\deg(\mathcal{C})-d+r)<d-r-1$.
Thus we have $2 < m < 3$ which contradicts that $m$ is an integer. 
Thus $m = 1$ and this implies that $\mathcal{C}$ must be a section.
\end{proof}
\begin{remark}\label{rem:delta1Cond}
If $\varphi:\mathbb{P}_{\mathbb{K}}^1\longrightarrow \mathcal{C}$ is birational onto its image (i.e., $\delta=1$), then the condition $d\geq \frac{2}{3}\bigl(\deg(\mathcal{C})+1\bigr)+r$ is equivalent to $d\geq 3r-2d'+2$.

On the other hand, the latter inequality $d\geq 3r-2d'+2$ forces $\delta=1$.
Indeed, since $\mathcal{C}\subset\mathbb{P}_{\mathbb{K}}^{d-r}$ is nondegenerate, we have $\deg(\mathcal{C})\geq d-r$. 
If $\delta\geq2$, then \Cref{cor:detectInnerProjections} gives $d-d'=\delta\deg(\mathcal{C})\geq2(d-r)$, and hence $d\leq2r-d'$.
Combining this with $d\geq3r-2d'+2$ yields $d'\geq r+2$ contradicting $d'\leq r$.
\end{remark}
As an additional result of the proof of \Cref{thm:surface-stable}, if $\mathcal{C}$ is a divisor on a smooth rational normal surface scroll $S(a,b)$, its Hilbert polynomial can be explicitly obtained.
\begin{remark}\label{rem:hilbertPolynomial}
Let $\mathcal{C}$ be a divisor on the smooth scroll $\mathcal{S} = S(a,b)$ of bidegree $(m,n)$ in $\mathrm{Pic}(\mathcal{S}) \cong \mathbb{Z}C_0 \oplus \mathbb{Z}F$. 
Then the Hilbert polynomial $\mathrm{HP}_{\mathcal{C}}(t)$ of $\mathcal{C}$ is given by
\[
\mathrm{HP}_{\mathcal{C}}(t) = \frac{d-d'}{\delta} t + \left( -mn + m + n + (b-a)\binom{m}{2} \right).
\]
\end{remark}
Note that \Cref{thm:minimalSectionScroll} gives an equivalent condition for $\rho(L)=r-1$.
Moreover, under the degree assumption $d\geq 3r-2d'+2$ (which forces $\varphi$ to be birational), we obtain a finer geometric characterization: the conditions $\rho(L)=r-1$ and $\max\{d_1,d_2\}<d$ are equivalent to $\mathcal C$ being smooth and lying on a smooth rational normal surface scroll.
\begin{corollary}\label{cor:surface-r}
Suppose that the rational curve $\mathcal{C}$ is not a rational normal curve.
Then $\rho(L)=r-1$ if and only if
\begin{enumerate}
    \item[\textnormal{(i)}] the morphism $\varphi:\mathbb{P}_{\mathbb{K}}^1\rightarrow\mathcal{C}$ defined by \eqref{equ:parametrizationCurve} is birational onto its image, and
    \item[\textnormal{(ii)}] $\mathcal{C}$ lies on a rational normal surface scroll in $\mathbb{P}_{\mathbb{K}}^{d-r}$, possibly a cone, as a section of its ruling.
\end{enumerate}
Moreover, assume that $d\geq 3r-2d'+2$.
Then $\rho(L)=r-1$ and $\max\{d_1,d_2\}<d$ if and only if
\begin{enumerate}
    \item[\textnormal{(iii)}] the curve $\mathcal{C}$ is smooth; and
    \item[\textnormal{(iv)}] $\mathcal{C}$ lies on a smooth rational normal surface scroll in $\mathbb{P}_{\mathbb{K}}^{d-r}$.
\end{enumerate}
\end{corollary}

\begin{proof}
The first statement follows immediately from \Cref{thm:minimalSectionScroll} with $m=2$.

Suppose that $d\geq 3r-2d'+2$.
By \Cref{rem:delta1Cond}, we have $\delta=1$.

Assume that $\rho(L)=r-1$ and $\max\{d_1,d_2\}<d$.
By \Cref{thm:minimalSectionScroll}, the curve $\mathcal{C}$ is contained as a section in the rational normal surface scroll $S(d-d_1,d-d_2)$.
Since both $d-d_1$ and $d-d_2$ are positive, this scroll is smooth.
Hence, by \Cref{thm:surface-stable}, the curve $\mathcal{C}$ is smooth.

Conversely, suppose that $\mathcal{C}$ is smooth and lies on a smooth rational normal surface scroll $\mathcal S=S(a,b)$.
Since $d\geq 3r-2d'+2$, the curve $\mathcal{C}$ is a section of the ruling of $\mathcal S$ by \Cref{thm:surface-stable}.
Together with $\delta=1$, the first statement gives $\rho(L)=r-1$.

It remains to show that $\max\{d_1,d_2\}<d$.
Suppose that $\max\{d_1,d_2\} = d$.
Without loss of generality, if $d_2=d$, then we have $d_1=r+1$ by \eqref{equ:totalDegree}.
By \Cref{thm:parametrizedCurve}, after a projective change of coordinates, $\varphi$ is given by
\[
\varphi([\alpha:\beta])=[g_1'(\alpha,\beta)\alpha^{d-r-1}:\cdots:g_1'(\alpha,\beta)\beta^{d-r-1}:g_2'(\alpha,\beta)].
\]
Since $\gcd(g_1',g_2')=1$, the inverse image of $[0:\cdots:0:1]$ is precisely $V(g_1')$.
Its length is $\deg(g_1')=r+1-d'\geq 2$, where the inequality follows from $d'\leq\rho(L)=r-1$.
On the other hand, since the morphism $\varphi$ is birational and $\mathcal{C}$ is smooth, $\varphi$ is an isomorphism.
Thus the inverse image of any point has length $1$ which is a contradiction.
Therefore, $\max\{d_1,d_2\}<d$ as desired.
\end{proof}
In particular, when $\dim(L) = 2$ and $d \geq 8$, a smooth curve $\mathcal{C} = \overline{\pi_{\mathcal{L}}(\widetilde{\mathcal{C}_d}\setminus \mathcal{L})}$ that is not a rational normal curve lies as a divisor on a smooth surface scroll if and only if $\rho(L) = 1$ by \Cref{cor:surface-r}. 
However, the following example demonstrates that the bound in \Cref{thm:surface-stable} is sharp.
\begin{example}\label{ex:surface-exception}
If $L= \langle S^5T^2, S^2T^5 \rangle$, then $\Ann(L)=(s^6,s^3t^3,t^6)$. 
The parametrization of the projected curve $\mathcal{C} = \overline{\pi_{\mathcal{L}}(\widetilde{\mathcal{C}_d}\setminus \mathcal{L})}$ is given by
\[
\{ [\alpha^7:\alpha^6\beta:\alpha^4\beta^3:\alpha^3\beta^4:\alpha \beta^6:\beta^7] \mid [\alpha:\beta] \in \mathbb{P}_{\mathbb{K}}^1 \} \subset\mathbb{P}_{\mathbb{K}}^5.
\]
Note that $p_a(\mathcal{C}) = 0$, but $d = 7 < 8 = 3r-2d'+2$.

Following \Cref{thm:minimalSectionScroll}, it is contained in $S(1,1,1)$.
However, it is contained in a surface scroll $S(2,2)$:
{\small
\begin{align*}
     \mathcal{C} &\cong \{[(\alpha^3)^2\alpha:(\alpha^3)(\beta^3)\alpha:(\beta^3)^2\alpha:(\alpha^3)^2\beta:(\alpha^3)(\beta^3)\beta:(\beta^3)^2\beta] \mid [\alpha:\beta] \in \mathbb{P}_{\mathbb{K}}^1\} \\
     & \subset  \{ [\gamma^2\eta:\gamma\delta \eta:\delta^2\eta:\gamma^2\zeta:\gamma\delta \zeta:\delta^2\zeta] \mid [\gamma:\delta],[\eta:\zeta] \in \mathbb{P}_{\mathbb{K}}^1\} = S(2,2).
\end{align*}
}
In fact, this is the case when $m=2+\frac{2r}{d-r-1}=2+\frac{4}{4}=3$. 
\end{example}
\begin{example}
If $L = \langle S^7T, ST^7\rangle$, then $\Ann(L)=(s^8, s^2t^2, t^8)$.
The parametrization of the projected curve $\mathcal{C} = \overline{\pi_{\mathcal{L}}(\widetilde{\mathcal{C}_8}\setminus \mathcal{L})}$ is given by
\[
\{ [\alpha^8:\alpha^6\beta^2:\alpha^5\beta^3:\alpha^4\beta^4:\alpha^3\beta^5:\alpha^2\beta^6:\beta^8] \mid [\alpha:\beta] \in \mathbb{P}_{\mathbb{K}}^1 \} \subset\mathbb{P}_{\mathbb{K}}^6.
\]

Following \Cref{thm:minimalSectionScroll}, it is contained in $S(0,0,4)$.
Note that $d = 3r-2d'+2 =8$, but $p_a(\mathcal{C}) = 2 >0$ and the scroll $S(0,0,4)$ is not smooth in this case.
In fact, the curve $\mathcal{C}$ is contained in the surface scroll $S(1,4)$, but it is not a section of its ruling. (In fact, $m=2$.)
More explicitly,
\begin{align*}
\mathcal{C} &\cong \bigl\{[(\alpha^2)^4:(\alpha^2)^3\beta^2:\cdots:(\beta^2)^4: \alpha^5\beta^3:\alpha^3\beta^5] \mid [\alpha:\beta]\in\mathbb P_{\mathbb K}^1 \bigr\} \\
&\subset \bigl\{[\gamma^4\eta:\gamma^3\delta\eta:\cdots:\delta^4\eta:\gamma\zeta:\delta\zeta]\mid [\gamma:\delta], [\eta:\zeta]\in\mathbb P_{\mathbb K}^1 \bigr\} = S(1,4).
\end{align*}
\end{example}

\subsection{Defining equations of rational curves on surface scrolls}
Curves of almost minimal degree that are not maximally del Pezzo are known to be linear projections of a rational normal curve in $\mathbb{P}^d$ from a point $p \in \mathbb{P}^d$ (cf.~\cite{MR2274517}). 
Such a curve $\mathcal{C}$ lies on a rational normal surface scroll, and its algebraic and geometric invariants including its Hilbert-Samuel polynomial, quadratic generation, and graded Betti numbers are determined by the smallest degree of minimal generators of $\operatorname{Ann}(F_p)$. 
Furthermore, it has recently been shown that the quadratic component of the defining ideal of $\mathcal{C}$ is spanned by quadratic forms of rank at most $4$ (for details, see \cite{MR4496651, MR5047032}).

In this section, we generalize these results to rational curves with $\rho(L) = r-1$.
We assume $d-r\ge2$ and order the two annihilator generators of degree at most $d$ so that $d_1\le d_2$. 
Thus $d_1+d_2=d+r+1$ and $r+1\le d_1\le d_2\le d$.
Suppose a rational curve $\mathcal{C} \subset \mathbb{P}_{\mathbb{K}}^{d-r}$ is projectively equivalent to the image of a rational normal curve $\widetilde{\mathcal{C}_d} \subset \mathbb{P}_{\mathbb{K}}^d$ under a linear projection from an $(r-1)$-dimensional linear subspace $\mathcal{L} \subset \mathbb{P}_{\mathbb{K}}^d$ with parametrization
{\small 
\begin{equation}\label{equ:paracurve}
  \mathcal{C} = \left\{ \left[ g_2'\alpha^{d-d_2} : \dots : g_2'\beta^{d-d_2} : g_1'\alpha^{d-d_1} : \dots : g_1'\beta^{d-d_1} \right] \mid [\alpha:\beta] \in \mathbb{P}_{\mathbb{K}}^1 \right\},  
\end{equation}
}
where $g_1', g_2' \in \mathbb{K}[\alpha,\beta]$ are binary forms of degree $d_1' = d_1 - d'$ and $d_2' = d_2 - d'$ with $d' = \deg(\gcd(g_1,g_2))$, respectively. 
Then, by \Cref{thm:minimalSectionScroll}, $\mathcal{C}$ is contained in the rational normal surface scroll $\mathcal{S} = S(d-d_2, d-d_1) \subset \mathbb{P}^{d-r}_{\mathbb{K}}$ with parametrization
{\small
\begin{equation}\label{equ:parasurfacescroll}
\mathcal{S} = \left\{ \left[ \gamma_1\alpha^{d-d_2} : \dots : \gamma_1\beta^{d-d_2} : \gamma_2\alpha^{d-d_1} : \dots : \gamma_2\beta^{d-d_1} \right] \mid [\alpha:\beta], [\gamma_1:\gamma_2] \in \mathbb{P}_{\mathbb{K}}^1 \right\}.
\end{equation}
}
Let $R(\mathcal{S})=\mathbb{K}[u,v,x,y]$ denote the total coordinate ring (i.e., Cox ring) of the standard projective-bundle model of $\mathcal{S}$ where $[u:v]$ and $[x:y]$ denote the base and local fiber coordinates, respectively. 
We equip $R(\mathcal{S})$ with the $\mathbb{Z}\times\mathbb{Z}$-grading given by
\[
\deg(u)=\deg(v)=(0,1), \quad \deg(x)=(1,d_2'), \quad \text{and} \quad \deg(y)=(1,d_1').
\]
Note that this grading differs from the usual convention, in which
\[
\deg(u)=\deg(v)=(0,1), \quad
\deg(x)=(1,d_2'-d_1'), \quad
\text{and} \quad \deg(y)=(1,0).
\]
In these coordinates, the section whose image is $\mathcal{C}$ is defined by the bihomogeneous equation $g_1'(u,v)x-g_2'(u,v)y=0$ of bidegree $(1,d_1'+d_2')=(1,d-2d'+r+1)$.
On the projective-bundle model, this section has divisor class $H+(r-d'+1)F$ where $H$ is the hyperplane class of $\mathcal{S}$ of bidegree $(1,d-d')$, and $F$ is the fiber class.

For integers $k \ge 0$ and $0 \le i \le k$, we define the bigraded vector subspaces $V_{k,i}$ and $W_{k,i}$ of the Cox ring $R(\mathcal{S})$ with the fiber term $x^{k-i}y^i$ fixed:
\[
V_{k,i} := x^{k-i}y^i \mathbb{K}[u,v]_{(k-i)(d-d_2)+i(d-d_1)-(r-d'+1)} \subset R(\mathcal{S})_{(k, k(d-d')-(r-d'+1))},
\]
\[
W_{k,i} := x^{k-i}y^i \mathbb{K}[u,v]_{(k-i)(d-d_2)+i(d-d_1)} \subset R(\mathcal{S})_{(k, k(d-d'))}.
\]
For $p \ge k$ and $i \le j \le i+p-k$, we say that $V_{p,j}$ is a \emph{multiple of $V_{k,i}$} if the natural multiplication map 
\[
V_{k,i} \otimes W_{p-k, j-i} \longrightarrow V_{p,j}
\]
is surjective. 
Since the multiplication map $\mathbb{K}[u,v]_a \otimes \mathbb{K}[u,v]_b \to \mathbb{K}[u,v]_{a+b}$ is surjective for all non-negative degrees $a, b \ge 0$, this map is surjective whenever both $V_{k,i}$ and $W_{p-k, j-i}$ are nonzero.

\begin{lemma}\label{lem:divisibility}
If $V_{k-1,i}$ is a nonzero vector space for some $k \geq 2$ and $0 \leq i \leq k-1$, then for any $p \geq k$ and any $j$ satisfying $i \leq j \leq i+p-k$, $V_{p-1,j}$ is a multiple of $V_{k-1,i}$.
\end{lemma}

\begin{proof}
By the transitivity of surjective multiplication maps, it suffices to establish the assertion for $p = k+1$ and $j \in \{i, i+1\}$. 
Since $V_{k-1,i}$ is non-trivial, the degree of its $\mathbb{K}[u,v]$-component must be non-negative, which implies that 
\[
(k-i-1)(d-d_2)+i(d-d_1)-(r-d'+1) \geq 0.
\]

First, suppose $p = k+1$ and $j = i$.
The coefficient degree of $V_{p-1,j}=V_{k,i}$ is
\begin{align*}
 & (k-i)(d-d_2)+i(d-d_1)-(r-d'+1) \\
= & \left[ (k-i-1)(d-d_2)+i(d-d_1)-(r-d'+1) \right] + (d-d_2).
\end{align*}
Since $d-d_2 \geq 0$, the multiplication map $V_{k-1,i} \otimes W_{1,0} \to V_{k,i}$ is surjective.

Next, suppose $p = k+1$ and $j = i+1$. 
The coefficient degree of $V_{p-1,j}=V_{k,i+1}$ is
\begin{align*}
 & (k-(i+1))(d-d_2)+(i+1)(d-d_1)-(r-d'+1) \\
= & \left[ (k-i-1)(d-d_2)+i(d-d_1)-(r-d'+1) \right] + (d-d_1).
\end{align*}
Since $d-d_1 \geq 0$, the multiplication map $V_{k-1,i} \otimes W_{1,1} \to V_{k,i+1}$ is surjective.

In both cases, $V_{p-1,j}$ is spanned by the image of $V_{k-1,i}$ under the respective multiplication maps, which completes the proof.
\end{proof}
By \Cref{lem:divisibility}, identifying the minimal generators of $I(\mathcal C)/I(\mathcal S)$ reduces to finding the nonzero components $V_{k-1,i}$ that are not obtained from lower-degree components by multiplication.

If $d_2=d$ (equivalently, $d_1=r+1$), the global sections defining the morphism are represented by $\{x,yu^{d-r-1},\dots,yv^{d-r-1}\}$.
In this case, $V_{k-1,i}=xV_{k-2,i}$ whenever $i<k-1$, so these components yield no new minimal generators.
Thus, we only need to determine the smallest integer $k$ for which $V_{k-1,k-1}$ is nonzero.
The relative ideal is then minimally generated by $(g_1'x-g_2'y)$ times a basis of this space.

Now, we prove that the minimal generators of $I(\mathcal{C})$ are given by differences of two products of linear forms, and determine the conditions for $I(\mathcal{C})$ to be generated by quadrics.
\begin{proof}[Proof of \Cref{thm:definingequs}]
The curve $\mathcal C$ in \eqref{equ:paracurve} lies as a divisor on the rational normal surface scroll $\mathcal S=S(d-d_2,d-d_1)$ given by \eqref{equ:parasurfacescroll}.
On the standard projective-bundle model of $\mathcal S$, the section mapping onto $\mathcal C$ is cut out by $g_1'x-g_2'y=0$ in $R(\mathcal S)$.
Identifying the coordinate ring of $\mathcal{S}$ with its hyperplane section subring $\mathbb{K}[xu^{d-d_2},\dots,xv^{d-d_2},yu^{d-d_1},\dots,yv^{d-d_1}]\subset R(\mathcal{S})$, the degree-$k$ component of the relative ideal $I(\mathcal{C})/I(\mathcal{S})$ is given by
\[
(I(\mathcal{C})/I(\mathcal{S}))_k \cong (g_1'x - g_2'y) \cdot \bigoplus_{i=0}^{k-1} V_{k-1,i}.
\]

\noindent\eqref{thm:definingequs_forms} 
Since $g_1',g_2'\in\mathbb K[u,v]$ factor into linear forms over the algebraically closed field $\mathbb K$, multiplying $g_1'x-g_2'y$ by any monomial basis element of $V_{k-1,i}$ yields a difference of two products of linear forms in $R(\mathcal S)$.
Thus, these expressions lift to differences of two products of $k$ linear forms in the homogeneous coordinates of $\mathbb P_{\mathbb K}^{d-r}$.
Together with the quadratic $2\times2$ minors defining $\mathcal S$, these lifts generate $I(\mathcal C)$, proving the claim.

\medskip
\noindent\eqref{thm:definingequs_quad} 
Quadratic generation holds if and only if $I(\mathcal{C})/I(\mathcal{S})$ is generated in degree $2$. 
If $d_1=r+1$, then $d-d_2=0$ and $V_{1,0}=0$. 
By \Cref{lem:divisibility}, the relative ideal is generated in degree $2$ exactly when $V_{1,1} \not= 0$.
The space $V_{1,1}$ is nonzero if and only if $(d-d_1)-(r-d'+1)=d-2r+d'-2\geq0$.
Thus quadratic generation is equivalent to $d\geq2r-d'+2$.

Now suppose $d_1>r+1$ (so both $d-d_1$ and $d-d_2$ are positive).
If $V_{1,0}$ and $V_{1,1}$ are both nonzero, \Cref{lem:divisibility} shows that all higher components are generated by them.
If $V_{1,0}=0$, then the first nonzero pure-$x$ component $V_{k-1,0}$ occurs for some $k \geq 3$. 
It cannot be generated by $V_{1,1}$ since their multiples all have positive $y$-exponent.
Therefore, quadratic generation is equivalent to $V_{1,0}\ne0$ since $d_1\le d_2$ implies that its dimension is at most that of $V_{1,1}$.
This condition is
\[
(d-d_2)-(r-d'+1)=d_1-2r+d'-2\ge0,
\]
or equivalently $d_1\ge2r-d'+2$.
\end{proof}
Although obtaining the complete numerical data for all minimal generators of $I(\mathcal{C})$ is complicated in general, it becomes highly tractable when $d_1 = r+1$ or $d \geq 3r-2d'+2$.
\begin{corollary}\label{cor:classification}
Under the same notation as above, the minimal generators of $I(\mathcal{C})$ are classified as follows:
\begin{enumerate}
    \item\label{cor:case1} If $d_1 = r+1$ (equivalently, $d_2=d$), then $I(\mathcal{C})$ is minimally generated by:
    \begin{enumerate}
        \item\label{cor:case1-1} $\binom{d-r}{2}-(r-d')$ quadratic forms if $d \ge 2r-d'+2$;
        \item\label{cor:case1-2} $\binom{d-r}{2}-(d-r-1)$ quadratic forms and $k(d-r-1)-d+d'+1$ forms of degree $k$ if $r+1 < d \leq 2r-d'+1$, where $k = \left\lceil \frac{r-d'+1}{d-r-1} \right\rceil+1 \ge 3$.
    \end{enumerate}

    \item\label{cor:case2} If $d \ge 3r-2d'+2$ and $d_1 > r+1$, then $I(\mathcal{C})$ is minimally generated by:
    \begin{enumerate}
        \item\label{cor:case2-1} $\binom{d-r}{2} - 2(r-d')$ quadratic forms if $d_1 \ge 2r-d'+2$;
        \item\label{cor:case2-2} $\binom{d-r}{2} - d_1' + 1$ quadratic forms and $k(d_1-r-1)-d_1'+1$ forms of degree $k$ if $d_1 \leq 2r-d'+1$, where $k = \left\lceil \frac{r-d'+1}{d_1-r-1} \right\rceil+1 \ge 3$.
    \end{enumerate}
\end{enumerate}
\end{corollary}
\begin{proof}
For each $k\in\mathbb N$, the exact sequence of vector spaces $0\to I(\mathcal S)_k\to I(\mathcal C)_k\to (I(\mathcal C)/I(\mathcal S))_k\to0$ gives $\dim I(\mathcal C)_k=\dim I(\mathcal S)_k+\dim(I(\mathcal C)/I(\mathcal S))_k$.
Also, $\dim I(\mathcal S)_2=\binom{d-r-1}{2}$.
For simplicity, set
\[
\delta_{k,i}:=(k-i)(d-d_2)+i(d-d_1)-(r-d'+1),
\]
so that $\dim_{\mathbb K}V_{k,i}=\max\{\delta_{k,i}+1,0\}$.

\medskip
\noindent\eqref{cor:case1}
When $d_1 = r+1$, the space $H^0(\mathcal S,\mathcal O_{\mathcal S}(1))$ is represented in the Cox ring $R(\mathcal S)$ by 
\[
\langle x, yu^{d-r-1}, \dots, yv^{d-r-1}\rangle.
\]
Therefore, the minimal degree of an additional generator is the first $k$ for which $V_{k-1,k-1} \not= 0$.

\eqref{cor:case1-1}
If $d \geq 2r-d'+2$, then $\delta_{(1,1)} = d-2r+d'-2 \geq 0$ which implies that no higher-degree generators are required.
The number of additional quadratic generators is $\delta_{(1,1)}+1 = d+d'-2r-1$. 

\eqref{cor:case1-2}
Suppose $r+1 < d < 2r-d'+2$.
Substituting $d_1 = r+1$ yields $\delta_{(k-1,k-1)} = (k-1)(d-r-1)-(r-d'+1)$. 
Thus, $\delta_{(k-1,k-1)} \geq 0$ if and only if $k \geq \left\lceil \frac{r-d'+1}{d-r-1} \right\rceil+1$ and the number of additional minimal generators is $\delta_{(k-1,k-1)}+1 = k(d-r-1)-d+d'+1$.

\medskip
\noindent\eqref{cor:case2}
From now on, we assume $d \geq 3r-2d'+2$ and $d_1>r+1$.

\eqref{cor:case2-1}
Assume $d_1 \geq 2r-d'+2$. 
The coefficient degree of $V_{1,i}$ is given by
\[
\delta_{(1,i)} = d-r-1-d_2'+i(d_2'-d_1').
\]
It is non-negative for both $i = 0$ and $i = 1$ because $d-r-1-d_2' = d_1-2r+d'-2 \geq 0$ and $d_1' \leq d_2'$. Therefore, we can construct a basis for $(g_1'x-g_2'y)(V_{1,0} \oplus V_{1,1})$ as follows:
\begin{align*}
    &\{u^{\delta_{(1,0)}-j} v^j g_1'(u,v) x^2 - u^{\delta_{(1,0)}-j} v^j g_2'(u,v) x y \;\mid\; 0 \leq j \leq \delta_{(1,0)}\} \\
    \cup \; &\{u^{\delta_{(1,1)}-j} v^j g_1'(u,v) x y - u^{\delta_{(1,1)}-j} v^j g_2'(u,v) y^2 \;\mid\; 0 \leq j \leq \delta_{(1,1)}\}.
\end{align*}
Consequently, $I(\mathcal{C})$ is minimally generated by $\dim I(\mathcal{S})_2 + (\delta_{(1,0)}+\delta_{(1,1)}+2) = \binom{d-r}{2} - 2(r-d')$ quadratic forms.

\eqref{cor:case2-2}
Assume $r+1 < d_1 < 2r-d'+2$. 
We have $\delta_{(1,0)} = d-r-d_2'-1 < 0$, but $\delta_{(1,1)} = d-r-d_1'-1 \geq 0$ since $d \geq 3r-2d'+2$ and $d_1' < 2r+2-2d'$. 
Therefore, by \Cref{lem:divisibility}, $V_{1,1}$ generates the components $V_{k-1,i}$ via multiplication with $W_{k-2,i-1}$ for all $1\leq i \leq k-1$.
Thus, the additional quadratic generators are given by
\[
\{u^{\delta_{(1,1)}-j} v^j g_1'(u,v) x y - u^{\delta_{(1,1)}-j} v^j g_2'(u,v) y^2 \;\mid\; 0 \leq j \leq \delta_{(1,1)}\},
\]
which contributes exactly $\delta_{(1,1)}+1 = d-r-d_1'$ elements to the scroll's quadrics.
Thus $I(\mathcal{C})_2$ is spanned by $\dim I(\mathcal{S})_2 + (\delta_{(1,1)}+1) = \binom{d-r}{2} - d_1' + 1$ quadratic forms.
To determine the remaining minimal generators, we find the smallest integer $k$ such that $V_{k-1,0}$ becomes non-trivial. 
This requires that
\[
(k-1)(d-d_2) - (r-d'+1) = (k-1)(d_1-r-1)-r+d'-1 \geq 0.
\]
This inequality yields the degree bound $k = \left\lceil \frac{r-d'+1}{d_1-r-1} \right\rceil + 1$, and the number of additional minimal generators is $\delta_{(k-1,0)} +1 = (k-1)(d_1-r-1)-r+d' = k(d_1-r-1)-d_1'+1$.
\end{proof}
More generally, suppose that $\mathcal C$ is contained as a section of the ruling of an $m$-fold rational normal scroll $\mathcal S=S(d-d_1,\dots,d-d_m)$ with $m\ge2$.
Following the convention above, write $R(\mathcal S)=\mathbb K[u,v,x_1,\dots,x_m]$ with $\deg(u)=\deg(v)=(0,1)$, $\deg(x_i)=(1,d_i-d')$, and hyperplane class $H=(1,d-d')$.

Let $M$ be a syzygy matrix of $(g_1,\dots,g_m)$, and set $J=([x_1\ \cdots\ x_m]M)$ and $I_{\mathrm{irr}}=(u,v)\cap(x_1,\dots,x_m)$.
The relative homogeneous ideal is obtained by saturating $J$ and taking its hyperplane-degree components:
\[
I(\mathcal C)/I(\mathcal S) \cong \bigoplus_{k\ge0}(J:I_{\mathrm{irr}}^\infty)_{kH}.
\]
In the surface case, $J=(g_1'x-g_2'y)$ is already saturated.
In higher codimension, saturation may introduce additional
minimal generators in the Cox ring which makes the ambient equations more difficult to determine.

\begin{example}
Suppose that $\mathcal{C} = \overline{\pi_{\mathcal{L}}(\widetilde{\mathcal{C}_{13}}\setminus \mathcal{L})}$ where $\mathcal{L} = \langle e_3, e_8 \rangle$.
Since $\operatorname{Ann}(S^{10}T^3, S^5T^8) = (s^{11}, s^6t^4, t^9)$, $\mathcal{C}$ is a codimension $2$ subvariety of the $3$-fold rational normal scroll $\mathcal{S} = S(2,3,4)$:
{\footnotesize
\[
\mathcal{S} = \left\{ [\gamma_1\alpha^2 : \cdots : \gamma_1\beta^2 : \gamma_2\alpha^3 : \cdots : \gamma_2\beta^3 : \gamma_3\alpha^4 : \cdots : \gamma_3\beta^4] \mid [\alpha:\beta] \in \mathbb{P}_{\mathbb{K}}^1, [\gamma_1:\gamma_2:\gamma_3] \in \mathbb{P}_{\mathbb{K}}^2 \right\}.
\]}
The local equations of $\mathcal{C}$ on $\mathcal{S}$ are given by
\[
v^4x_1 - u^5x_2 = 0 \quad \text{and} \quad v^5x_2 - u^6x_3 = 0,
\]
which are obtained by computing the minimal syzygies of the annihilator generators $[u^{11}\;\;\; u^6v^4\;\;\; v^9]$.

Even though the local ideal is generated by only two simple equations, computing the saturation with respect to the irrelevant ideal $I_{\mathrm{irr}} = (u,v) \cap (x_1,x_2,x_3)$ produces additional generators in the Cox ring.
More explicitly, using the software system Macaulay2 \cite{M2}, we can find the following additional generators through the saturation:
\[
u x_1 x_3 - v x_2^2, \quad u^4 x_2^3 - v^3 x_1^2 x_3, \quad u^3 x_2^5 - v^2 x_1^3 x_3^2, \quad u^2 x_2^7 - v x_1^4 x_3^3, \quad u x_2^9 - x_1^5 x_3^4.
\]
\end{example}

\section*{Acknowledgement}
The authors thank Euisung Park for valuable discussions that helped clarify the significance of the results presented in this paper.
The first and second named authors were supported by the National Research Foundation of Korea (NRF) grant funded by the Korean government (MSIT) (No.~RS-2024-00414849).
The second named author was also supported by an NRF grant funded by MSIT (No.~RS-2026-25482052).
The third named author was supported by the Sungshin Women’s University Research Grant of 2025(H20250019).

\subsection*{AI Statement}
The conceptual framework and principal mathematical arguments presented in this work were developed by the authors without the use of artificial intelligence.
During the authors' review of the manuscript, ChatGPT (OpenAI; GPT-6 Astra) was consulted to assist with checking minor technical details and relevant references, and improving English expression and the presentation of the manuscript.
The authors independently evaluated all AI-assisted suggestions and take full responsibility for the accuracy, integrity, and final content of the work.

\bibliographystyle{alpha}
\bibliography{reference.bib}

\end{document}